\documentclass[11pt,a4paper]{article}
\usepackage[T1]{fontenc}
\usepackage[utf8]{inputenc}
\usepackage{amsmath}
\usepackage{authblk}
\usepackage{graphicx}
\usepackage{indentfirst}
\usepackage{float}
\usepackage{amsfonts}
\usepackage{cite}
\usepackage[colorlinks,citecolor=blue]{hyperref}
\usepackage{xcolor}
\usepackage{geometry}
\usepackage{lineno}
\newtheorem{theorem}{Theorem}
\newtheorem{lemma}[theorem]{Lemma}
\date{}
\title{\bf{Entropy stable non-oscillatory fluxes: An optimized wedding of entropy conservative flux with non-oscillatory flux}}
\author[1,2]{\bf{Ritesh Kumar Dubey}\thanks{riteshkd@srmist.edu.in, ritesh@blockappsai.com}\thanks{Author acknowledges SERB India for fund through project EMR/2016/000394 to support authors research visit to Blockapps AI, Bangalore where initial work for this paper got carried out. Author also acknowledges SERB India for project CRG/2022/002659 which enabled author to complete this manuscript and further developments.}}
\affil[1]{\textit{Research Institute,Department of Mathematics,SRM Institute of Science and Technology, Chennai, India\\
}}
\affil[2]{\textit{Blockapps AI, Bangalore, India}}

\begin{document}

\maketitle

\begin{abstract}
This work frames the problem of constructing  non-oscillatory entropy stable fluxes as a least square optimization problem. A flux sign stability condition is defined for a pair of entropy conservative flux ($F^*$) and a non-oscillatory flux ($F^s$). This novel approach paves a way to construct non-oscillatory entropy stable flux ($\hat{F}$) as a simple combination of $(F^*$ and $F^s)$ which inherently optimize the numerical diffusion in the entropy stable flux ($\hat{F}$) such that it reduces to the underlying non-oscillatory flux ($F^s$) in the flux sign stable region.	This robust approach is (i) agnostic to the choice of flux pair $(F^*,F^s)$, (ii) does not require the computation of costly dissipation operator and high order reconstruction of scaled entropy variable to construct the diffusion term. Various non-oscillatory entropy stable fluxes are constructed and exhaustive computational results for standard test problems are given which show that fully discrete schemes using these entropy stable fluxes do not exhibit nonphysical spurious oscillations in approximating the discontinuities compared to the non-oscillatory schemes using underlying fluxes ($F^s$) only. Moreover, these entropy stable schemes maintain the formal order of accuracy of the lower order flux in the pair.

\end{abstract}
{\bf \textit{Keywords:}} Hyperbolic conservation laws, Entropy stability, Maximum principle, High order non-oscillatory schemes, Sign stability property, Least square optimization.\\
{\bf \textit{AMS subject classifications:}} 65M06, 35L65
\section{Introduction}
One of the topics of paramount interest in the area of numerical solution of partial differential equations is to design schemes for non-linear hyperbolic systems of conservation laws. These equations are natural models for many physical problems in areas of engineering and scientific studies like  gas dynamics, fluid flow problems, astrophysical flow etc. Generic one dimensional system of conservation laws can be written in the form,
\begin{eqnarray}\label{conservation_laws}
\dfrac{\partial \mathbf{u}}{\partial t} +\dfrac{\partial \mathbf{f}(\mathbf{u})}{\partial x }&=&0,\;\;	\forall (x,\;t) \in   \mathbb{R} \times \mathbb{R}^+,\\
\mathbf{u}(x,0) &=& \mathbf{u}_0(x), \forall x \in \mathbf{R}, 
\end{eqnarray}
where vector $\mathbf{f}(\mathbf{u}):\mathbb{R}^m \rightarrow \mathbb{R}^m $  is a smooth flux function of conserved vector quantity  $\mathbf{u}(x,~t):\mathbb{R}\times \mathbb{R}^+\rightarrow \mathbb{R}^m$. 
In order to numerically approximate the solution of \eqref{conservation_laws}, one discretizes the spatial domain into intervals $I_{i}=[x_{i-\frac{1}{2}},x_{i+\frac{1}{2}}]$ using Cartesian mesh $\{x_i\}_{i\in \mathbb{Z}}$ with mesh size $\Delta x= x_{i+1}-x_{i}$ such that $x_{i}= (x_{i+\frac{1}{2}} +x_{i-\frac{1}{2}})/2$. The semi-discrete finite difference (or finite volume) conservative scheme is given by \begin{equation}
		\label{semi_scheme}
		\frac{d}{dt}\mathbf{u}_{i}(t)=-\frac{1}{\Delta x_{i}}\left({\mathbf{F}}_{i+\frac{1}{2}}-{\mathbf{F}}_{i-\frac{1}{2}}\right),
	\end{equation}
	where $\mathbf{u}_{i}(t)=\mathbf{u}(x_{i}, t)$ and numerical flux function  ${\mathbf{F}}_{i+\frac{1}{2}}$ satisfies the  consistency condition, that is ${\mathbf{F}}(\mathbf{u},\mathbf{u},...,\mathbf{u})=\mathbf{f}(\mathbf{u})$. 
The computation of the updated solution of \eqref{conservation_laws} requires a fully discrete scheme which really obtained using a suitable time discretization in \eqref{semi_scheme} e.g., \cite{gottlieb1998total}. Note that, until mentioned explicitly, solution using scheme \eqref{semi_scheme} should be understood as solution by the fully discrete version of \eqref{semi_scheme} depending on underlying time marching process as discussed in section \ref{sec6}.
\par It is well known that \eqref{conservation_laws} admits discontinuities like shock and rarefaction in its solution which pose further problems like non-unique weak solutions,  crisp and efficient numerical resolution of these discontinuities. Some of the well known classical schemes converge to physically incorrect weak  solution \cite{leveque2002finite}. On one hand low order schemes yield smeared approximation of discontinuities whereas on other hand uniformly high order schemes exhibit spurious oscillations around such discontinuities. Therefore, it is required to device robust and efficient high order schemes which crisply capture these discontinuities without introducing spurious oscillations. Moreover, one core requirement is that such numerical schemes should yield solution which converges to viscosity solution of \eqref{conservation_laws}\cite{evans1998partial}. 

 \par Notably the viscosity solution shows two intrinsic properties (i) it uniquely satisfies entropy stability inequality \cite{Tadmor2003} and (ii) in scalar case, it satisfies a strict maximum principle \cite{Zhang2011}. These characteristics of the viscosity solution paved the way for devising modern shock capturing schemes. 
 There is a vast literature on numerical schemes which are designed based on above characteristics of physically acceptable solution. We mention a few foundational as well recent significant work which fall inline with the theme of this paper.
 \par  Non-oscillatory schemes which are designed following the maximum principle are monotone schemes \cite{crandall1979monotone}, high resolution total variation diminishing (TVD) schemes \cite{Harten1984, sweby1984high,toro2000,Zhang2011, ZhangTVD,Kumar2007483}. Schemes that weakly satisfy maximum principle are arbitrary high order accurate Essentially non-oscillatory (ENO) \cite{SHU1988439}, weighted ENO  \cite{Shu1997} schemes and their modified improved versions in \cite{jiang1996efficient, serna2004power, castro2011high, fan2014new, kim2016modified,rathan2018modified, henrick2005mapped, biswas2018accuracy} etc. Though such maximum principle satisfying schemes ensure non-oscillatory approximation of discontinuous solution, they do not guarantee convergence to viscosity solution even in scalar case. It is important to note that there are schemes which are entropy stable as well as strictly satisfy maximum principle e.g., monotone and TVD schemes \cite{Sanders1983,Ismail2009, DUBEY2018,Chen2011,Chen2017}, however they are at most first order accurate at solution extrema \cite{osher1984high, osher1986very, tadmor1988convenient}. Note that, conservative three-point second-order accurate methods cannot satisfy a local entropy inequality  \cite{maria}.
 
\par Entropy stable fluxes for hyperbolic conservation laws can be obtained by adding numerical diffusion into entropy conservative fluxes \cite{Tadmor1987}.   Some of the entropy stable schemes designed following this approach are \cite{Tadmor1987, Ismail2009,Fjordholm2012,Zakerzadeh19052015,Cheng2019, DUBEY2018,BbRk,Fisher2013,Cheng2016,DUAN2021110136, Winters201672}. Note that, entropy stable schemes guarantee for convergence to viscosity solution but they do not ensure for non-occurrence of oscillations in the vicinity of discontinuities particularly higher order schemes due to the absence of suitable numerical diffusion \cite{BbRk}. Note that deducing the suitable amount of this additive numerical diffusion such that resulting entropy stable scheme maintains formal high order of accuracy and yet yields non-oscillatory approximate solution is non-trivial and has been a topic of wide interest and research. The existing approach for constructing high order non-oscillatory entropy stable flux using high order numerical diffusion operator demands for (i) sign stability of the reconstructed scaled entropy variable which is very restrictive requirement (ii) explicit computation of dissipation operator which is often costly and problem dependent.

\par The aim of this work is to construct robust and efficient arbitrary high order non-oscillatory entropy stable scheme for hyperbolic conservation laws \eqref{conservation_laws}. The idea is to design non-oscillatory entropy stable fluxes such that additive numerical diffusion (i) is implicitly defined (ii) is robust in the sense that it does not depend on sign property of the high order reconstructions of scaled entropy variable (iii) is optimized in such a way that the non-oscillatory property of entropy stable flux is governed by some established underlying non-oscillatory flux. 

\par The rest of the paper is as follows: For motivation and completeness of the presentation a brief introduction and challenges in the construction of {\em high order non-oscillatory entropy stable fluxes} are highlighted in section \ref{sec2}. In \ref{sec3}, a quick introduction on non-oscillatory fluxes like TVD and ENO/WENO is given. The main contribution of this work is given in section \ref{sec4} where task of constructing non-oscillatory entropy stable flux is posed as least square optimization problem and required numerical diffusion is deduced using first order optimality condition. In section \ref{sec5}, a flux sign stability lemma for entropy stability of any flux is given and used to retrospectively analyze some well known entropy stable fluxes. In section \ref{sec6}, numerical results by various constructed entropy stable schemes are given and compared. Conclusions are given in section \ref{sec7}.    
 \section{Entropy stable scheme} \label{sec2}
 This section briefly present theory and development of entropy stable schemes; and comprehensive details can be found in \cite{Tadmor2003, Tadmor2016}. We assume that system \eqref{conservation_laws} is equipped with entropy pair $\left(\eta(\mathbf{u}), q(\mathbf{u})\right)$ which symmetrizes the system \eqref{conservation_laws}. Then it can be shown that a physical weak solution of \eqref{conservation_laws} uniquely satisfies following entropy inequality
\begin{equation}\label{Entropy_condition}
\frac{\partial}{\partial t} \eta(\mathbf{u})+\frac{\partial}{\partial x}q(\mathbf{u}) \leq 0,
\end{equation}
where entropy function $\eta \equiv \eta(\mathbf{u}):\mathbb{R}^m\rightarrow \mathbb{R}$ is convex  and entropy flux function $q(\mathbf{u}): \mathbb{R}^m\rightarrow \mathbb{R} $ satisfies the following compatibility relation with entropy variable $\mathbf{v}=\eta_\mathbf{u}(\mathbf{u})$
\begin{equation}\label{sym_cond}
\mathbf{v}^T \mathbf{f_u}=q_\mathbf{u}^T.
\end{equation}
By virtue of the convexity of $\eta$, the mapping $\mathbf{v}\rightarrow \mathbf{u}$ is one-to-one and change of variable $\mathbf{u} =  \mathbf{u(v)}$ transforms the system \eqref{conservation_laws} into an equivalent symmetric form 
\begin{equation}\label{sym_form}
\frac{\partial}{\partial t} \mathbf{u(v)}+\frac{\partial}{\partial x}\mathbf{g}(\mathbf{v})=0,\; \mathbf{g(v)} =\mathbf{f(u(v))}.
\end{equation} 
The system \eqref{sym_form} is symmetric in the sense that the Jocobians of temporal and spatial fluxes of variable $\mathbf{v}$ i.e., $\mathbf{u(v)}$ and $\mathbf{g(v)}$ respectively satisfy
\begin{equation}\label{jacob_mat}
H(\mathbf{v})= \mathbf{u_v(v)}= H^T(\mathbf{v})>0, \;B(\mathbf{v}) =\mathbf{g_v(v)}=B^T(\mathbf{v}).
\end{equation}
Moreover, under \eqref{sym_cond}, the Hessian of an entropy function symmetrizes the system  \eqref{conservation_laws} as follows \cite{Tadmor2003}
\begin{equation}
\eta_{\mathbf{uu}}A = \left[\partial \eta_{\mathbf{uu}}A\right]^T,\; A = \dfrac{\partial f(\mathbf{u}) }{\partial \mathbf{u}} .
\end{equation}
The semi-discrete scheme \eqref{semi_scheme} is said to be entropy conservative if the computed solution satisfies the following discrete entropy criterion,
\begin{equation}\label{disEC}
\frac{d}{dt}\eta(\mathbf{u}_{i}(t))+\frac{1}{\Delta x_{i}}\left({q}_{i+\frac{1}{2}}-{q}_{i-\frac{1}{2}}\right)=0,
\end{equation}
where discrete entropy flux function ${q}_{i+\frac{1}{2}}={q}(\mathbf{u}_{i-l+1},...,\mathbf{u}_{i+l}), l \in \mathbb{Z}^+$ is consistent with entropy flux function $q$ such that
\begin{equation}
{q}(\mathbf{u},\mathbf{u},...,\mathbf{u})=q(\mathbf{u}).
\end{equation}
Define the notations for jump and average of a discrete quantity ${z}_{i},\; i\in \mathbb{Z}$ by  $\left[\left[z\right]\right]_{i}:=z_{i+\frac{1}{2}}-z_{i-\frac{1}{2}}$ and  $\bar{z}_{i+\frac{1}{2}} := \frac{1}{2}(z_{i+1}+z_{i})$ respectively. Then a foundational approach  to construct entropy conservative numerical flux in \cite{Tadmor1987} is as follows. 
\begin{theorem} \label{theorm1987}
Let a consistent numerical flux  $ {\mathbf{F}}_{i+\frac{1}{2}}={\mathbf{F}}_{i+\frac{1}{2}}^*$ satisfies 
\begin{equation}\label{eq_ec}
\left[\left[\mathbf{v}\right]\right]_{i+\frac{1}{2}} {\mathbf{F}}_{i+\frac{1}{2}}^*=\left[\left[\psi\right]\right]_{i+\frac{1}{2}},
\end{equation}
where $\psi $ is entropy potential given by 
\begin{equation} \label{defpsi}
\psi(\mathbf{v})=\mathbf{v}\cdot \mathbf{g}(\mathbf{v})-q(\mathbf{u}(\mathbf{v})).
\end{equation}
Then the semi-discrete scheme \eqref{semi_scheme} with numerical flux $\mathbf{F}^*$ is second order accurate and entropy conservative i.e., solution computed by scheme satisfies the discrete entropy equality \eqref{disEC} with numerical entropy flux 
\begin{equation}\label{qstar}
q_{i+\frac{1}{2}} \equiv {q}_{i+\frac{1}{2}}^*(\mathbf{u}_{i}, \mathbf{u}_{i+1}) = \bar{\mathbf{v}}^{T}_{i+\frac{1}{2}} \mathbf{F}^{*}_{i+\frac{1}{2}} - \bar{\psi}_{i+\frac{1}{2}}.
\end{equation}
\end{theorem}
The entropy conservative flux ${\mathbf{F}}^*$ obtained from \eqref{eq_ec} for scalar problem is unique, though suffers from non uniqueness in case of the system of conservation laws. However, this non-uniqueness does not severely impact the development of entropy conservative fluxes. Various second as well arbitrary high order entropy conservative fluxes can be found in \cite{Chandrashekar2013,Ismail2009,Winters201672,Tadmor2003, LeFloch2002, Fjordholm2012}. A detailed comparison of various entropy conservative fluxes for Euler equation is given in \cite{Hendrik2017}. 
\subsection{Entropy Stable fluxes and numerical diffusion} The entropy of the solution of hyperbolic conservation law \eqref{conservation_laws} remains conserved for smooth solution and dissipates only across shocks \cite{BbRk}. Therefore, schemes using entropy conservative fluxes work well for smooth solutions however allow non-physical oscillations near shock due to the absence of numerical diffusion which ensures that computed solution exhibits entropy dissipation. A standard procedure to obtain such dissipative fluxes  is to add numerical diffusion explicitly in to the entropy conservative flux \cite{Tadmor1987, Tadmor2003}. We state the result therein as follows,
\begin{theorem}\label{ES_Thm}
	Let $\mathbf{D}_{i+\frac{1}{2}}\in\mathbf{R}^{m\times m}$ be symmetric positive semi-definite i.e., $\mathbf{D}_{i+\frac{1}{2}}\geq 0$  and $\mathbf{F}^{*}:\mathbf{R}^m\rightarrow \mathbf{R}^m$ be entropy conservative flux, then the solution by semi-discrete scheme \eqref{semi_scheme} using numerical flux of the form 
	\begin{equation}\label{ES_Flux}
	\hat{\mathbf{F}}_{i+\frac{1}{2}}=\mathbf{F}_{i+\frac{1}{2}}^*-\frac{1}{2} \mathbf{D}_{i+\frac{1}{2}}\left[\left[\mathbf{v}\right]\right]_{i+\frac{1}{2}},\,
	\end{equation}
satisfies the following discrete analog of the entropy inequality \eqref{Entropy_condition},
\begin{equation}\label{disECInq}
\frac{d}{dt}\eta(\mathbf{u}_{i}(t))+\frac{1}{\Delta x_{i}}\left[\hat{q}_{i+\frac{1}{2}}-\hat{q}_{i-\frac{1}{2}}\right]\leq 0,
\end{equation}
where numerical entropy flux $\hat{q}_{i+\frac{1}{2}} = q^{*}_{i+\frac{1}{2}} -\frac{1}{2}\bar{\mathbf{v}}^{T}_{i+\frac{1}{2}} \mathbf{D}_{i + \frac{1}{2}}\left[\left[\mathbf{v}\right]\right]_{i+\frac{1}{2}}$. 
\end{theorem}
The defined flux $\hat{\mathbf{F}}_{i+\frac{1}{2}}$ is called entropy stable flux and the  term $\frac{1}{2}\mathbf{D}_{i+\frac{1}{2}}\left[\left[\mathbf{v}\right]\right]_{i+\frac{1}{2}}$ in \eqref{ES_Flux} represents numerical diffusion term. The numerical diffusion term in \eqref{ES_Flux} is constituted by a symmetric positive semi-definite dissipation operator $\mathbf{D}_{i+\frac{1}{2}}$ and the jump in discrete entropy variable $\mathbf{v}_i$.  Note that,"A conservative scheme which contain more numerical viscosity than that present in the entropy conservative one is also entropy stable" \cite{Tadmor1987}. 
Thus, constructing flux $\hat{\mathbf{F}}_{i+\frac{1}{2}}$ which achieve entropy stabilty can trivially be done using any symmetric $\mathbf{D}_{i+\frac{1}{2}}\geq 0$. However, achieving {\em high order accuracy and non-oscillatory property} in the entropy stable fluxes is non trivial. In particular:
\subsubsection{Quest for suitable dissipation operator} 
The dissipation of the numerical entropy determines the discontinuity capturing ability of the underlying scheme. Thus just the positiveness property of dissipation operator $\mathbf{D}_{i+\frac{1}{2}}$ alone is not enough to ensure for truly oscillation free entropy stable scheme. Further, excessive dissipation causes smeared approximation of the discontinuity whereas in the absence of sufficient diffusion, spurious oscillations may occur in the vicinity of discontinuities.  Thus the problem is to determine a suitable diffusion operator $\mathbf{D}_{i+\frac{1}{2}}$ which enables the scheme to yield a crisp and non-oscillatory resolution of discontinuities. 
\par Let $\mathbf{A}=\partial_{\mathbf{u}}\mathbf{f(u)}$ be the Jacobian matrix with complete set of independent eigen vectors such that $\mathbf{A}= \mathbf{R}\mathbf{\Lambda} \mathbf{R}^{-1}$, where $\Lambda$ is a non-negative diagonal matrix depending on eigenvalues of Jacobian $\mathbf{A}$ and $\mathbf{R}$ is matrix of associated eigenvectors. Some examples of diffusion operator used in \cite{Fjordholm2012} are $\mathbf{D}_{i+\frac{1}{2}} = \tilde{\mathbf{R}}_{i+\frac{1}{2}}\tilde{\mathbf{\Lambda}}_{i+\frac{1}{2}} \tilde{\mathbf{R}}^T_{i+\frac{1}{2}}$, where $\mathbf{\tilde{R}}$ is a scaling of $\mathbf{R}$ such that $\mathbf{u_v}= \mathbf{\tilde{R}\tilde{R}}^T$ and  $\mathbf{A}= \mathbf{\tilde{R}}\mathbf{\Lambda} \mathbf{\tilde{R}^{T}}$.
The subscript $i+\frac{1}{2}$ of $\mathbf{\tilde{R}}$ and $\mathbf{\tilde{\Lambda}}$ denotes an average state at cell interface $x_{i+\frac{1}{2}}$. Two of the choices used for $\mathbf{\Lambda}$ therein are Roe type dissipation operator $\mathbf{\Lambda} = diag(|\lambda_1|, |\lambda_2|, \dots |\lambda_m|)$ and Rusanov type dissipation operator  $\mathbf{\Lambda} = \max\left(|\lambda_1|, |\lambda_2|, \dots |\lambda_m|\right)\mathbf{I}_{m\times m}$. 
In \cite{Ismail2009}, entropy consistent diffusion operator {\it EC1} is proposed using a second and third-order differential terms for Euler equation and it is given as,
\begin{equation}\label{ec1}
\mathbf{\Lambda}^{EC1} = |\mathbf{\Lambda}| + \frac{1}{6}|\left[\left[\mathbf{\Lambda}_{u\pm a}\right]\right]|,
\end{equation}
 where $\left[\left[\mathbf{\Lambda}\right]\right] = diag\left(\left[\left[ u-a \right]\right],0, \left[\left[u+a\right]\right]\right)$. However, extension of \eqref{ec1} for general hyperbolic systems is not clear and its non-oscillatory nature depends on the choice of entropy function \cite{DUBEY2018}. A total variation diminishing condition is deduced on  $\mathbf{D}$ such that the resulting entropy stable flux \eqref{ES_Flux} ensures for the complete removal of spurious oscillations in \cite{DUBEY2018}. For systems it is defined as,
 \begin{equation}\label{estvd}
 \mathbf{D}_{i+\frac{1}{2}} \geq \max\left(\mathbf{\tilde{R}}|\mathbf{\Lambda}|\mathbf{\tilde{R}}^T + |Q^*|, 0\right),
 \end{equation}
 where the entropy viscosity matrix $Q^*$ in terms of Jacobian $B(\mathbf{v})$ in \eqref{jacob_mat} is given by,
 \begin{equation}\label{ev_matrix}
 Q^{*}= \int_{-\frac{1}{2}}^{\frac{1}{2}}2\xi\,B\left(\mathbf{v}_{i + \frac{1}{2}}(\xi)\right),\; \,\mathbf{v}_{i+\frac{1}{2}}=\bar{\mathbf{v}}_{i+\frac{1}{2}} + \xi \left[\left[\mathbf{v}\right]\right]_{i+\frac{1}{2}}.
 \end{equation}
 Clearly, to obtain diffusion operator in \eqref{estvd}, one needs to explicitly calculate entropy viscosity which often may not be efficiently computable. Apart from above, the diffusion operator $\mathbf{D}$ in the context of Shallow water hydrodynamics equation is given by using Cholesky decomposition of matrix $\dfrac{\partial \mathbf{u}}{\partial \mathbf{v}}$ in \cite{DUAN2021110136}. Thus it can be observed that for general system \eqref{conservation_laws}, defining and computing a suitable $\mathbf{D}$ is a complicate need to design non-oscillatory entropy stable flux. 
        
\subsubsection{Need for sign stable high order reconstruction}  Note that irrespective of order of accuracy of flux $\mathbf{F}^{*}$, entropy stable schemes with numerical diffusion term $\frac{1}{2}\mathbf{D}_{i+\frac{1}{2}}\left[\left[\mathbf{v}\right]\right]_{i+\frac{1}{2}}$ are only first order accurate as the jump $\left[\left[\mathbf{v}\right]\right]_{i+\frac{1}{2}}$ in entropy variable across cell interface $x_{i+\frac{1}{2}}$ is of order $\Delta x$. Thus in order to achieve high order entropy stable flux, a high order reconstruction of the jump in the entropy variable $\mathbf{v}$ is needed as shown in TECNO schemes \cite{Fjordholm2012}. More precisely, let $\mathbf{v}_{i}^{\pm}$ be $(2k-1)^{th}$ order reconstruction of entropy variable $\mathbf{v}$, then $(2k-1)^{th}$ order entropy stable flux is obtained by adding  $(2k-1)^{th}$ order diffusion term to $2k^{th}$ order entropy conservative flux in the form
\begin{equation}\label{tecno_f}
\hat{\mathbf{F}}_{i+\frac{1}{2}}=\mathbf{F}_{i+\frac{1}{2}}^*-\frac{1}{2} \tilde{\mathbf{R}}_{i+\frac{1}{2}}\tilde{\mathbf{\Lambda}}_{i+\frac{1}{2}}\left<\left<\mathbf{\tilde{w}}\right>\right>_{i+\frac{1}{2}},\,
\end{equation}

where $\left<\left<\mathbf{\tilde{w}}\right>\right>_{i+\frac{1}{2}} = \mathbf{\tilde{w}}_{i+1}^+ - \mathbf{\tilde{w}}_{i}^-$ is the jump in the reconstructed values of scaled entropy variable $\mathbf{w}^{\pm}_i= \tilde{\mathbf{R}}^T_{i+\frac{1}{2}}\mathbf{v}_i$ defined as  $ \mathbf{\tilde{w}}^{\pm}_i = \tilde{\mathbf{R}}^T_{i\pm\frac{1}{2}}\mathbf{v}^{\pm}_i$. Note that,
\begin{subequations}\label{wv}
\begin{equation}
	\left<\left<\mathbf{\tilde{w}}\right>\right>_{i+\frac{1}{2}} = \mathbf{\tilde{w}}_{i+1}^+ - \mathbf{\tilde{w}}_{i}^- = \mathbf{\tilde{R}}_{i+\frac{1}{2}}^T\mathbf{v}_{i+1}^- - \mathbf{\tilde{R}}_{i+\frac{1}{2}}^T\mathbf{v}_{i}^+ = \mathbf{\tilde{R}}_{i+\frac{1}{2}}^T\left<\left<\mathbf{v}\right>\right>_{i+\frac{1}{2}},\label{tildew}
	\end{equation}
\begin{equation}
	\left<\left<\mathbf{w}\right>\right>_{i+\frac{1}{2}} = \mathbf{{w}}_{i+1}^+ - \mathbf{w}_{i}^- = \mathbf{\tilde{R}}_{i+\frac{1}{2}}^T\mathbf{v}_{i+1} - \mathbf{\tilde{R}}_{i+\frac{1}{2}}^T\mathbf{v}_{i} \;=\; \mathbf{\tilde{R}}_{i+\frac{1}{2}}^T\left[\left[\mathbf{v}\right]\right]_{i+\frac{1}{2}}.\label{wonly}
\end{equation}	
\end{subequations}
It is shown in \cite{Fjordholm2012} that flux \eqref{tecno_f} is entropy stable provided reconstruction of scaled entropy variable  satisfies the following component wise sign property at each interface $x_{i+\frac{1}{2}}$,
\begin{equation}\label{Tencosignprop}
sign\left(\left<\left<\mathbf{\tilde{w}}\right>\right>_{i+\frac{1}{2}}\right) = sign\left(\left<\left<\mathbf{w}\right>\right>_{i+\frac{1}{2}}\right),
\end{equation} 
For each component $l$, \eqref{Tencosignprop} is defined as
\begin{equation}\label{componantSS}
	\begin{array}{cc}
		\left<\left<w^l\right>\right>_{i+\frac{1}{2}}>0, & \mbox{then}\; \left<\left<\tilde{w}^l\right>\right>_{i+\frac{1}{2}}\geq 0,\\
		\left<\left<w^l\right>\right>_{i+\frac{1}{2}}<0, & \mbox{then}\; \left<\left<\tilde{w}^l\right>\right>_{i+\frac{1}{2}}\leq 0,\\
		\left<\left<w^l\right>\right>_{i+\frac{1}{2}}=0, & \mbox{then}\; \left<\left<\tilde{w}^l\right>\right>_{i+\frac{1}{2}}= 0.\\
	\end{array}
\end{equation}
In scalar case, \eqref{Tencosignprop} reduces to 
\begin{equation}\label{scalarsignprop}
	sign\left(\left<\left<{v}\right>\right>_{i+\frac{1}{2}}\right)= sign\left(\left[\left[{v}\right]\right]_{i+\frac{1}{2}}\right).
\end{equation}
It is needed to state here that, scaling in \eqref{wv} involve expansive matrix-vector multiplication between scaled matrix $\mathbf{\tilde{R}}$ and entropy vector $\mathbf{v}$ at each cell interface.

Further, note that TeCNO framework \cite{Fjordholm2012} for constructing high order non-oscillatory schemes demands for sign stability property \eqref{Tencosignprop} in the high order reconstruction and therefore, only ENO reconstruction could be used therein. Authors in \cite{BbRk} modified the TECNO framework which demands sign stablity of scaled entropy variable only accross the  locally significantly jumps and hence can work with other high order reconstructions e.g., third order WENO and high order TVD reconstruction. In \cite{DUAN2021110136}, high-order accurate well-balanced semi-discrete entropy stable schemes are developed for shallow water magneto hydrodynamics by adding diffusion using a switch function proposed in \cite{BbRk}. The construction of suitable dissipation term therein again based on the WENO reconstruction of the scaled entropy variables.  Recently in \cite{LIU2019104266}, the third order WENO reconstruction of scaled entropy variable is proposed to construct entropy stable schemes for shallow water equations.
 
\section{Non-Oscillatory Schemes}\label{sec3}
Here we present a quick review of class of prevailing non-oscillatory schemes. Note that, oscillatory approximations for discontinuous solution of \eqref{conservation_laws}  can not be considered
as admissible solution even for scalar case since it violates the following global maximum
principle MP satisfied by the its physically admissible solution. In scalar case MP is given as,
\begin{equation} \min_{x}(u_{0}(x)) \leq u(x,t) \leq \max_{x}(u_{0}(x)),
\forall (x,\;t)\in \mathbb{R}\times\, \mathbb{R}^{+}. 
\label{mp} 
\end{equation}
We call a scheme non-oscillatory in the sense that the discrete solution satisfies either the maximum principle \eqref{mp} or the number of extrema of the discrete solution is not increasing in time as defined in \cite{HartenOsher}. We also take the liberty to call essentially non-oscillatory (ENO) schemes as non oscillatory though they do not follow the maximum principle. The rationale for doing so lies in the fact that the proposed framework can lead to Monotone, TVD or essentially non-oscillatory entropy stable scheme depending on the choice of entropy conservative flux and such non-oscillatory flux (See sub-sections \ref{sec5a} and \ref{sec5b}).
As mentioned earlier, examples of non-oscillatory schemes which Maximum principle \eqref{mp} satisfying schemes are monotone schemes \cite{crandall1979monotone}, total variation diminishing (TVD) schemes \cite{Harten1984,sweby1984high, ZhangTVD,Goodman}. Apart from these schemes, other non-oscillatory schemes which do not strictly follow maximum principle \cite{Zhang2011} but practically give excellent non-oscillatory numerical results are uniformly non-oscillatory scheme \cite{HartenOsher}, essentially non-oscillatory (ENO) and weighted ENO schemes \cite{SHU1988439,Shu1997, Shu2009} and references therein. Among these non-oscillatory schemes, ENO and WENO schemes are very attractive as they preserve formal higher order of accuracy unlike monotone and TVD schemes and are well developed now \cite{jiang1996efficient, henrick2005mapped,borges2008improved,castro2011high,ha2013improved,rathan2018modified,biswas2018accuracy,Sabana}. 
\subsection{High order non-oscillatory fluxes} 
The fluxes corresponding to ENO/WENO schemes rely on high order reconstruction/interpolation of the conserved quantity at cell interface $x_{i+\frac{1}{2}}$. Let  $v(x),\; x\in \Sigma \subset \mathbb{R}$ be a piece-wise continuous function and the domain $\Sigma$  be partitioned with the grids $\{x_i\}$, $i\in\mathbb{Z}$, and the point values be given by $v_{i}=v(x_{i})$. Then the  $k$-th order ENO interpolation procedure in an arbitrary interval $I_{i}:=[x_{i-\frac{1}{2}}, x_{i+\frac{1}{2}}]$ utilizes $(2k-1)$ grid point stencil and consists of two steps. The first step  chooses smoothest stencil among $k$ consecutive points $S_{pref}=\{x_{i-r},\,\dots,\,x_{i},\,\dots,\, x_{i-r+k-1}\}$, where $r\in\{0,1,\dots, \left(k-1\right)\}$. Then, a unique $\left(k-1\right)$-th degree polynomial $p_i(x)$ passing through such $S_{pref}$ is used to interpolate conserved variable at cell interfaces $x_{i\pm\frac{1}{2}}$ as, 
\begin{eqnarray*}
	v_{i+\frac{1}{2}}^{-}=p_{i}(x_{i+\frac{1}{2}}),
	\\v_{i-\frac{1}{2}}^{+}=p_{i}(x_{i-\frac{1}{2}}).
\end{eqnarray*}
 Compared to ENO, with WENO reconstruction one can achieve improved $(2k-1)^{th}$ order accuracy using same $(2k-1)$ point stencil for smooth data. The idea in WENO reconstruction is to consider a convex combination of reconstructed values at cell interfaces $x_{i+\frac{1}{2}}$ using all $(k-1)^{th}$ order unique polynomials $p_i^r(x)$ which reconstruct function $v(x)$ over sub stencil $S_r(i)=\{x_{i+r-k+1},......,x_{i+r}\},\;\;r=0,....k-1$\cite{Shu2009}. Such $(2k-1)^{th}$ order accurate reconstructed values are given by 
\begin{equation}\label{WENO_Re}
v_{i+\frac{1}{2}}^{-}=\sum_{r=0}^{k-1} \omega_r p^{r}_i(x_{i+\frac{1}{2}}),\;
v_{i-\frac{1}{2}}^{+}=\sum_{r=0}^{k-1} \tilde{\omega}_r p^{r}_i(x_{i-\frac{1}{2}}),
	\end{equation}
where non-linear weights $w_r, \tilde{w}_r$ are given by
\begin{subequations}
	\begin{equation}
	\omega_r=\frac{\alpha_r}{\sum_{p=0}^k\alpha_p},\, \, \tilde{\omega}_r=\frac{\tilde{\alpha}_r}{\sum_{p=0}^k\tilde{\alpha}_p},
	\end{equation}
	with 
	\begin{equation}\label{beta_j}
	\alpha_r=\frac{\gamma_{r}}{(\epsilon+\beta_r)^2},\,\,\tilde{\alpha}_r=\frac{\tilde{\gamma}_{r}}{(\epsilon+\beta_r)^2}.
	\end{equation}
\end{subequations}
The constants $\gamma_{r}$ and $\tilde{\gamma}_r$ are such that 
\begin{equation*}
\sum_{r=0}^{k-1} \gamma_r p^{r}_i(x_{i+\frac{1}{2}})-v(x_{i+\frac{1}{2}})=O(h^{2k-1}),
\end{equation*}
and
\begin{equation*}
\sum_{r=0}^{k-1} \tilde{\gamma}_r p^{r}_i(x_{i-\frac{1}{2}})-v(x_{i-\frac{1}{2}})=O(h^{2k-1}).
\end{equation*}
The parameters $\beta_r$'s in \eqref{beta_j} measure the smoothness and are given by\begin{equation}\label{beta}
\beta_r=\sum_{l=1}^{k}\int_{x_{i-\frac{1}{2}}}^{x_{i+\frac{1}{2}}}{\Delta x}^{2l-1}\left(\frac{d^l}{dx^l}p^{j}_i(x)\right)^2 dx, \,(j=0,1,...,k-1).
\end{equation} 
A good detail on ENO/WENO reconstruction and interpolation procedure can be found in \cite{Shu1997, Fjordholm2013}. In scalar case, the  ENO/WENO flux $F^{s}_{i+\frac{1}{2}}$ in the finite difference semi-discrete conservative approximation \eqref{semi_scheme} is reconstructed by point values $f(u_i)$. More precisely, use $f^+(u_i)$ to reconstruct the positive cell interface numerical flux  $F^{+}_{i+\frac{1}{2}}= v^{-}_{i+\frac{1}{2}}$ and use $f^{-}(u_i)$ to reconstruct the negative cell interface numerical flux  $F^{-}_{i+\frac{1}{2}}= v^{+}_{i+\frac{1}{2}}$. The positive and negative flux $f^{\pm}$ satisfy $f(u) = f(u)^+ +f(u)^-$  which can be obtained by any suitable flux splitting such that $\frac{f(u)^+}{du} \geq 0,\; \frac{f(u)^-}{du}\leq 0$ \cite{Shu1997}. For example, following simple Lax-Friedrichs splitting,
\begin{equation}\label{split_lxf}
f^{\pm}(u) = \frac{1}{2}(f(u) \pm \sigma u),
\end{equation} 
where $\sigma = \displaystyle \max_u|f'(u)|$ over the relevant range of $u$.
Finally, the non-oscillatory ENO/WENO numerical flux can be obtained by $F^{s}_{i+\frac{1}{2}} = F^{+}_{i+\frac{1}{2}}+ F^{-}_{i+\frac{1}{2}}$.
We conclude this section by stating that despite of active efforts to develop high order entropy stable schemes, entropy stable fluxes are never designed to construct schemes which explicitly mimic non-oscillatory properties of schemes obtained by established aforementioned non-oscillatory fluxes. 
 
\section{Least square optimization Problem}\label{sec4}

Throughout in the paper, vector functions $\hat{\mathbf{F}}$ and $\mathbf{F}^*$ denote generic numerical entropy stable flux and entropy conservative flux respectively whereas $\mathbf{F}^s$ denotes non-oscillatory numerical flux function of any stable scheme like monotone, TVD or ENO, WENO schemes etc. More precisely, let $\mathbf{F}^s$ be any $n^{th}$ order consistent non-oscillatory numerical flux for \eqref{semi_scheme} of a non-oscillatory finite difference scheme such that 
\begin{equation}\label{2kEC}
\frac{1}{\Delta x}(\mathbf{F}_{i+\frac{1}{2}}^s - \mathbf{F}_{i-\frac{1}{2}}^s )= \left. \frac{\partial \mathbf{f(\mathbf{u})}}{\partial x}\right|_{x=x_i} + \mathcal{O}(\Delta x^{n}).
\end{equation}
Also, let $\mathbf{F}^*$ be a $m^{th}$ order consistent entropy conservative flux for \eqref{ES_Flux} such that 
\begin{equation}\label{2kFw}
\frac{1}{\Delta x}(\mathbf{F}_{i+\frac{1}{2}}^* - \mathbf{F}_{i-\frac{1}{2}}^* )= \left.\frac{\partial \mathbf{f(\mathbf{u})}}{\partial x}\right|_{x=x_i} + \mathcal{O}(\Delta x^{m}).
\end{equation} 
Further, the scheme with flux $\mathbf{F}^*$ satisfies the discrete entropy equality \eqref{disEC} with some appropriate numerical entropy function $q$. For example, in $m=2p^{th}$, $p \in \mathbb{N}$ order entropy conservative schemes \cite{LeFloch2002} the numerical entropy function is
\begin{equation}\label{2p_qstar}
q_{i+\frac{1}{2}} =  \sum_{k=1}^p \alpha_k^p\sum_{l=1}^{k-1}q^*_{i+\frac{1}{2}}(\mathbf{u}_{i-l}, \mathbf{u}_{i-l+k}),
\end{equation}
where $q^*$ is second order numerical entropy flux function \eqref{qstar} and $\alpha_k^p$'s solve the linear system 
\begin{equation}
\sum_{k=1}^p k\alpha_k^p =1, \; \sum_{i=1}^{p} i^{2j -1}\alpha_k^p =0, (j= 2, 3,\dots p).
\end{equation}

\begin{lemma}\label{lem-new}
	Let $\mathbf{u}, \mathbf{v} \in \mathbb{R}^m$ and $sign(\mathbf{u})= sign(\mathbf{v})$ i.e, $sign({u_i})= sign({v_i})\; \forall i = 1,2, \dots m$, then there exists a is symmetric positive definite diagonal scaling matrix $\mathbf{S} =diag(s_{11}, s_{22}, \dots, s_{mm}) $ where $s_{ii}\geq0, \forall i = 1,2,\dots,m$ such that $\mathbf{Su} = \mathbf{v}$. 
\end{lemma}
\par Note that for the problem \eqref{conservation_laws}, these numerical fluxes are discrete vector functions such as $\hat{\mathbf{F}}_{i+\frac{1}{2}}: \mathbb{R}^m \rightarrow \mathbb{R}^m  $, 
$\mathbf{F}^*_{i+\frac{1}{2}}: \mathbb{R}^m \rightarrow \mathbb{R}^m  $,
and $\mathbf{F}^s_{i+\frac{1}{2}}: \mathbb{R}^m \rightarrow \mathbb{R}^m  $. Denote by $\mathbb{P}$ the set of symmetric positive definite matrices $\textbf{P}\in \mathbb{R}^{m\times m}$ s.t. $ \textbf{x}^T\textbf{P} x\geq0, \forall \textbf{x} \in \mathbb{R}^m$ and $M_{i,j}$ represents $(i,j)^{th}$ element of a matrix $\mathbf{M}$. Further, denote by $\nabla_{\textbf{D}}$ the matrix differentiation operator with respect to matrix $\mathbf{D}_{i+\frac{1}{2}}$ for example defined either in $\omega$ or $\alpha$-derivative sense \cite{MAGNUS20102200}.  

In this setting, using \eqref{ES_Flux} the problem of constructing entropy stable flux $\hat{\mathbf{F}}_{i+\frac{1}{2}}$ which can imitate non-oscillatory property of $\mathbf{F}^s_{i+\frac{1}{2}}$, be posed as following modified least square problem i.e.,\\
$$\mathcal{P}1: \mbox{Given}\; \left[\left[\mathbf{v}\right]\right]_{i+\frac{1}{2}} \in  \mathbb{R}^m, \mbox{find matrix} \;\mathbf{D}_{i+\frac{1}{2}} \in \mathbb{P} \subset \mathbb{R}^{m\times m} \; s.t.\; \frac{1}{2}\mathbf{D}_{i+\frac{1}{2}}\left[\left[\mathbf{v}\right]\right]_{i+\frac{1}{2}} = \mathbf{F}_{i+\frac{1}{2}}^*-\mathbf{F}_{i+\frac{1}{2}}^s.$$
Using \eqref{ES_Flux}, $\mathcal{P}1$ can be written as: find matrix
$\mathbf{D}_{i+\frac{1}{2}} \in \mathbb{P} \subset \mathbb{R}^{m\times m} \; s.t.\; \hat{\mathbf{F}}_{i+\frac{1}{2}} = \mathbf{F}_{i+\frac{1}{2}}^s.$ This least square problem can be reformulated as following equivalent minimization problem 
\begin{equation}\label{LSprob1}
\mathcal{P}2: \displaystyle \min_{\mathbf{D} \in \mathbb{P}} \frac{1}{2}\|\hat{\mathbf{F}}_{i+\frac{1}{2}} - \mathbf{F}^s_{i+\frac{1}{2}}\|^2 
\end{equation}
The problem of deducing such minimizing $\mathbf{D}$ can equivalently be written as optimization problem \begin{equation}	
\mathcal{P}3: \;\;\mathbf{D} = \displaystyle  \underset{\mathbf{D}_{i+\frac{1}{2}}\in \mathbb{P}}{\arg\min} J(\mathbf{D}_{i+\frac{1}{2}})\; \mbox{where}\;  J(\mathbf{D}_{i+\frac{1}{2}}) = \frac{1}{2}\left(\mathbf{F}^*_{i+\frac{1}{2}} -\frac{1}{2} \mathbf{D}_{i+\frac{1}{2}}\left[\left[\mathbf{v}\right]\right]_{i+\frac{1}{2}} -\mathbf{F}^s_{i+\frac{1}{2}}\right)^2,
\end{equation}
where $J(\mathbf{D}_{i+\frac{1}{2}})$ is the cost or penalty function. Thanks to convexity of $J(\mathbf{D}_{i+\frac{1}{2}})$ a unique minimizer $\mathbf{D}_{i+\frac{1}{2}}$ exists and can be determined by using the first order optimality condition. Differentiating $J(\mathbf{D})$ w.r.t. matrix $\mathbf{D}$ in $\omega$-derivative sense, we have,
\begin{equation}\label{gradJD}
\nabla_{\mathbf{D}_{i+\frac{1}{2}}} J(\mathbf{D})= \left(\mathbf{F}^*_{i+\frac{1}{2}} -\frac{1}{2} \mathbf{D}_{i+\frac{1}{2}}\left[\left[\mathbf{v}\right]\right]_{i+\frac{1}{2}} -\mathbf{F}^s_{i+\frac{1}{2}}\right)\left(-\frac{1}{2}\nabla_{\mathbf{D}_{i+\frac{1}{2}}}(\mathbf{D}_{i+\frac{1}{2}}\left[\left[\mathbf{v}\right]\right]_{i+\frac{1}{2}})\right)
\end{equation}
Consider the first order sufficient optimality condition i.e.,  $\nabla_{\mathbf{D}_{i+\frac{1}{2}}} J(\mathbf{D}) =0$ to deduce the minimizing optimizer $\mathbf{D}^o$. Following two cases arise
	\begin{itemize}
		\item Let $\left[\left[\mathbf{v}\right]\right]_{i+\frac{1}{2}}=\mathbf{0}$, then the tensor $\nabla_{\mathbf{D}_{i+\frac{1}{2}}}\left(\mathbf{D}_{i+\frac{1}{2}}\left[\left[\mathbf{v}\right]\right]_{i+\frac{1}{2}}\right) \equiv \mathbf{0}$  identically. This case corresponds a non-shock solution region where entropy remains conservative and entropy stability is trivially achieved by $\mathbf{F}^*_{i+\frac{1}{2}}$ (irrespective of the choice of $\mathbf{D}$ in \eqref{ES_Flux}). In fact, in this case from minimization problem it follows that, $\mathbf{F}^*_{i+\frac{1}{2}} = \mathbf{F}^s_{i+\frac{1}{2}}$. It also explains the well behaved numerical approximation by entropy conservative schemes in the case of smooth solution. 
		\item Let $\left[\left[\mathbf{v}\right]\right]_{i+\frac{1}{2}}\neq \mathbf{0}$, then the tensor $\nabla_{\mathbf{D}_{i+\frac{1}{2}}}\left(\mathbf{D}_{i+\frac{1}{2}}\left[\left[\mathbf{v}\right]\right]_{i+\frac{1}{2}}\right)\not \equiv \mathbf{0}$ identically. Then \eqref{gradJD} imposes the following condition on minimizing optimizer $\mathbf{D}^o$ 
		 \begin{equation}\label{num_diff1}
		\frac{1}{2} \mathbf{D}^o_{i+\frac{1}{2}}\left[\left[\mathbf{v}\right]\right]_{i+\frac{1}{2}} = \mathbf{F}^*_{i+\frac{1}{2}}-\mathbf{F}^s_{i+\frac{1}{2}}.
		\end{equation}
	\end{itemize}
Solving \eqref{num_diff1} to obtain minimizing  $\mathbf{D}^o$ in system case is non-trivial. However, the good thing to be observed from \eqref{ES_Flux} is that, to construct entropy stable flux $\hat{\mathbf{F}}_{i+\frac{1}{2}}$ the explicit computation of minimizing $\mathbf{D}^o$ is not needed and diffusion term $\frac{1}{2} \mathbf{D}^o_{i+\frac{1}{2}}[\mathbf{v}]_{i+\frac{1}{2}}$ can be directly replaced by right hand side expression of \eqref{num_diff1}. In addition, for entropy stability of the flux $\hat{\mathbf{F}}$,  the positivity condition $\mathbf{D}^o\geq 0$ must hold (refer Theorem \ref{ES_Thm}). Note from lemma \ref{lem-new} and under the following  {\em flux sign stability property}. 
\begin{equation}\label{signprop}
	sign\left(\mathbf{F}^*_{i+\frac{1}{2}}-\mathbf{F}^s_{i+\frac{1}{2}}\right) = sign\left(\left[\left[\mathbf{v}\right]\right]_{i+\frac{1}{2}}\right).
\end{equation}
The unique minimizing symmetric positive diagonal matrix $\mathbf{D}^o\geq 0$ can be defined as
\begin{equation*}
	d_{ii} = \left\{\begin{array}{cc}\displaystyle	
	\frac{F^*_i -F^s_i}{v_{i+1}-v_{i}}, & \,\text{if}\, sign(F^*_i -F^s_i)=sign(v_{i+1}-v_{i})\neq 0\\ &\\
	 d_{ij} =0 & \text{else}  \end{array}\right.
\end{equation*}
$\forall i,j = 1,2,\dots, m$. Further, using property \eqref{signprop}, the {\em optimized non-oscillatory numerical diffusion} for constructing  flux \eqref{ES_Flux} is defined as 
\begin{equation}\label{num_diff2}
	\frac{1}{2} \mathbf{D}^o_{i+\frac{1}{2}}\left[\left[\mathbf{v}\right]\right]_{i+\frac{1}{2}}=
	\left\{\begin{array}{cc}
		\mathbf{F}^*_{i+\frac{1}{2}}-\mathbf{F}^s_{i+\frac{1}{2}} &\; \mbox{if}\; sign\left(\mathbf{F}^*_{i+\frac{1}{2}}-\mathbf{F}^s_{i+\frac{1}{2}}\right) = sign\left(\left[\left[\mathbf{v}\right]\right]_{i+\frac{1}{2}}\right), \\
		\mathbf{0} &\; \mbox{else}.
	\end{array}\right.
\end{equation}
Alternatively, the entropy stable flux $\hat{\mathbf{F}}_{i+\frac{1}{2}}$ can directly be expressed as the following convex combination
\begin{equation}\label{eq_chi}
\hat{\mathbf{F}}_{i+\frac{1}{2}} = \mathbf{F}^*_{i+\frac{1}{2}} - \chi_{i+\frac{1}{2}} \left(\mathbf{F}^*_{i+\frac{1}{2}}-\mathbf{F}^s_{i+\frac{1}{2}}\right),
\end{equation}
where 
\begin{equation}\label{chi_eq} 
\chi_{i+\frac{1}{2}} = \left\{\begin{array}{cc}
1 & \mbox{if}\;sign\left(\mathbf{F}^*_{i+\frac{1}{2}}-\mathbf{F}^s_{i+\frac{1}{2}}\right) = sign\left(\left[\left[\mathbf{v}\right]\right]_{i+\frac{1}{2}}\right),  \\
0 & \;\mbox{else}.
\end{array}\right.
\end{equation}

\begin{theorem}
The numerical scheme with flux \eqref{ES_Flux} and numerical diffusion \eqref{num_diff2} is entropy stable and its solution satisfies 
\begin{equation}\label{disESInq}
\frac{d}{dt}\eta(\mathbf{u}_{i}(t))+\frac{1}{\Delta x_{i}}\left[\hat{q}_{i+\frac{1}{2}}-\hat{q}_{i-\frac{1}{2}}\right]\leq 0,
\end{equation}	
where discrete entropy function
$\displaystyle \hat{q}_{i+\frac{1}{2}} = q^*_{i+\frac{1}{2}} -\frac{1}{2}\bar{\mathbf{v}}^{T}_{i+\frac{1}{2}} \mathbf{D}^o_{i + \frac{1}{2}}\left[\left[\mathbf{v}\right]\right]_{i+\frac{1}{2}}$ and $q^*_{i+\frac{1}{2}}$ is the numerical entropy flux associated with the entropy conservative flux $\mathbf{F}^*_{i+\frac{1}{2}}$.\\
{\bf Proof:} It follows from \eqref{num_diff2} that
\begin{equation}
\left[\left[\mathbf{v}\right]\right]_{i+\frac{1}{2}}^T \mathbf{D}^o_{i+\frac{1}{2}}\left[\left[\mathbf{v}\right]\right]_{i+\frac{1}{2}}=
\left\{\begin{array}{cc}
\left[\left[\mathbf{v}\right]\right]^T_{i+\frac{1}{2}} \left(\mathbf{F}^*_{i+\frac{1}{2}}-\mathbf{F}^s_{i+\frac{1}{2}}\right) &\; \mbox{if}\; sign\left(\left[\left[\mathbf{v}\right]\right]_{i+\frac{1}{2}}\right) = sign\left(\mathbf{F}^*_{i+\frac{1}{2}}-\mathbf{F}^s_{i+\frac{1}{2}}\right) \\
\mathbf{[v]}^T_{i+\frac{1}{2}}\cdot \mathbf{0} &\; \mbox{else}
\end{array}\right\}  \geq \mathbf{0}.
\end{equation}
Now on multiplying the scheme \eqref{semi_scheme} with $\mathbf{v}^T$ and following the steps as in \cite{Tadmor1987} we get 
\begin{eqnarray}
\frac{d}{dt}\eta(\mathbf{u}_{i}(t)) +\frac{1}{\Delta x_{i}}\left[\hat{q}_{i+\frac{1}{2}}-\hat{q}_{i-\frac{1}{2}}\right]&=&
-\frac{1}{4\Delta x}\left(
\left[\left[\mathbf{v}\right]\right]_{i+\frac{1}{2}}^T \mathbf{D}^o_{i+\frac{1}{2}}\left[\left[\mathbf{v}\right]\right]_{i+\frac{1}{2}} + \left[\left[\mathbf{v}\right]\right]_{i-\frac{1}{2}}^T \mathbf{D}^o_{i-\frac{1}{2}}\left[\left[\mathbf{v}\right]\right]_{i-\frac{1}{2}}\right)\nonumber \\
&\leq & \mathbf{0}.\nonumber
\end{eqnarray}
\hfill $\square$
\end{theorem}
\section{Non-oscillatory Entropy stable flux}\label{sec5}
In scalar case, on applying Mean value theorem on entropy variable ${v({u})}$ we get ,  
\begin{equation}
\displaystyle \left[\left[{v}\right]\right]_{i+\frac{1}{2}}={v}_{{u}}(\xi)\left[\left[{u}\right]\right]_{i+\frac{1}{2}},
\end{equation}
where $\xi\in [{u}_i,\,{u}_{i+1}]$. Since $\eta({u})$ is convex, ${v_u}({\xi})=\eta_{{uu}}({\xi}) \geq 0$ and it follows
\begin{equation}\label{sign_uisv}
sign\left(\left[\left[{v}\right]\right]_{i+\frac{1}{2}}\right)=sign\left(\left[\left[{u}\right]\right]_{i+\frac{1}{2}}\right).
\end{equation}
For systems, the component-wise non-oscillatory entropy stable flux can be defined using \eqref{eq_chi} in the following simple form
\begin{equation}\label{finalESflux}
\mathbf{\hat{F}}_{i+\frac{1}{2}}=
\left\{\begin{array}{cc}
\mathbf{F}^s_{i+\frac{1}{2}} &\; \mbox{if}\; sign\left(\left[\left[\mathbf{v}\right]\right]_{i+\frac{1}{2}}\right) = sign\left(\mathbf{F}^*_{i+\frac{1}{2}}-\mathbf{F}^s_{i+\frac{1}{2}}\right), \\
\mathbf{F}^*_{i+\frac{1}{2}} &\; \mbox{else}.
\end{array}\right.
\end{equation}
The following  trivially follows for system \eqref{conservation_laws} which is equipped with symmetrizing entropy pair $\left(\eta(\mathbf{u}), q(\mathbf{u})\right)$  and provides an entropy stability condition for a given flux using {\em flux sign stability property} \eqref{signprop}.
\begin{lemma}\label{lemm2}
Let $\mathbf{F}$ be any consistent numerical flux for system \eqref{conservation_laws} and if there exists a consistent entropy conservative flux $\mathbf{F}^*$ for \eqref{conservation_laws} such that
\begin{equation}\label{F_isES}
\mathbf{F}_{i+\frac{1}{2}} =\mathbf{F}^*_{i+\frac{1}{2}} - \frac{1}{2}\mathbf{\tilde{D}}_{i+\frac{1}{2}}\left[\left[\mathbf{v}\right]\right]_{i+\frac{1}{2}}, 
\end{equation}  
where $\mathbf{\tilde{D}}_{i+\frac{1}{2}}$ is any unknown diagonal scaling matrix, then $\mathbf{F}$ is entropy stable provided it satisfies the following component-wise flux sign stability property 
\begin{equation}\label{signpropu}
sign\left(\mathbf{F}^*_{i+\frac{1}{2}}-\mathbf{F}_{i+\frac{1}{2}}\right) = sign\left(\left[\left[\mathbf{v}\right]\right]_{i+\frac{1}{2}}\right).
\end{equation}
{\bf Proof:} Rewrite \eqref{F_isES} as
\begin{equation}
\mathbf{F}^*_{i+\frac{1}{2}} -\mathbf{F}_{i+\frac{1}{2}} = \frac{1}{2}\mathbf{\tilde{D}}_{i+\frac{1}{2}}\left[\left[\mathbf{v}\right]\right]_{i+\frac{1}{2}}, 
\end{equation}
Clearly, under the flux sign stability conditions \eqref{signpropu}, scaling matrix  $\mathbf{\tilde{D}}_{i+\frac{1}{2}}\geq 0$. Thus from Theorem \ref{ES_Thm} the flux \eqref{F_isES} is entropy stable.
\hfill $\square$
\end{lemma}

\pagebreak

\par{\bf Remark:} {\it Lemma \ref{lemm2} shows that the second branch in the switching function $\chi$ in \eqref{chi_eq} gets chosen except in case if non-oscillatory flux $F^s$ itself is entropy stable.}
\par In the following, Lemma \ref{lemm2} is tested for some well known non-oscillatory entropy stable fluxes including high order TeCNO fluxes. It is shown that under flux sign stability property \eqref{signpropu} they satisfy the entropy stability condition.   
\subsection{First order Entropy stable flux}\label{sec5a}
We remark that the first order entropy stable \textcolor{blue}{TVD} fluxes in \cite{Ismail2009, DUBEY2018}, are of the form \eqref{F_isES} and under \eqref{sign_uisv} naturally satisfy flux sign stability condition \eqref{signpropu} for the dissipation operators characterizes by \eqref{ec1} or \eqref{estvd}. \textcolor{black}{We further consider the following flux of a generic three point conservative scheme,
\begin{equation}\label{llf_u}
	\mathbf{F}_{i+\frac{1}{2}} = \bar{\mathbf{f}}_{i+\frac{1}{2}} -\frac{1}{2} \alpha_{i+\frac{1}{2}}\left[\left[\mathbf{u}\right]\right]_{i+\frac{1}{2}},
\end{equation}
where \begin{equation}\label{alpha_cond}
	\alpha_{i+\frac{1}{2}} \geq  \max_{\lambda,|\xi|\leq 1/2}\left|\lambda\left[\mathbf{A(u(v}_{i+\frac{1}{2}}(\xi)))\right]\right|.
\end{equation}
On writing the second order entropy conservative flux \eqref{eq_ec} in viscosity form \cite{Tadmor2003},
\begin{equation}\label{visc_form}
	\mathbf{F}^{*}_{i+\frac{1}{2}} = \bar{\mathbf{f}}_{i+\frac{1}{2}} - \dfrac{1}{2}Q^*_{i+\frac{1}{2}} \left[\left[\mathbf{v}\right]\right]_{i+\frac{1}{2}},
\end{equation} 
where $Q^*$ is defined in \eqref{ev_matrix}.  
Note from \eqref{jacob_mat} that 
\begin{equation}
	\left[\left[\mathbf{u}\right]\right]_{i+\frac{1}{2}}= \int_{\xi =-1/2}^{1/2} \dfrac{d}{d\xi}\mathbf{u(v}(\xi))_{i+\frac{1}{2}} d\xi =
	\int_{\xi =-1/2}^{1/2}H(v_{i+\frac{1}{2}}(\xi))d\xi. \left[\left[\mathbf{v}\right]\right]_{i+\frac{1}{2}}.
\end{equation}
Thus using $\alpha_{i+\frac{1}{2}}\left[\left[\mathbf{u}\right]\right]_{i+\frac{1}{2}} = Q_{i+\frac{1}{2}} \left[\left[\mathbf{v}\right]\right]_{i+\frac{1}{2}}$, flux \eqref{llf_u} can be written in terms of entropy variable,
\begin{equation}\label{llf_v}
	\mathbf{F} = \bar{\mathbf{f}}_{i+\frac{1}{2}} -\frac{1}{2} Q_{i+\frac{1}{2}}\left[\left[v\right]\right]_{i+\frac{1}{2}},
\end{equation}
where
\begin{equation}
	Q_{i+\frac{1}{2}} = \int_{\xi =-1/2}^{1/2}H(v_{i+\frac{1}{2}}(\xi))d\xi .
\end{equation}
From \eqref{visc_form} and \eqref{llf_v} we have 
$$\mathbf{F}^*_{i+\frac{1}{2}}-\mathbf{F}_{i+\frac{1}{2}} = \frac{1}{2}(Q-Q^*)\left[\left[\mathbf{v}\right]\right]_{i+\frac{1}{2}}.$$ 
Let $\mathbf{F}^*$ and $\mathbf{F}$ satisfies the flux sign stability property \eqref{signpropu} then from lemma \ref{lem-new}, there exists symmetric positive definite diagonal matrix $\mathbf{S}= Q_{i+\frac{1}{2}} - Q^{*}_{i+\frac{1}{2}}\geq 0$. Therefore, the flux $\mathbf{F}_{i+\frac{1}{2}}$ in \eqref{llf_u} is entropy stable.
We note that in \cite{Tadmor2003} using viscosity comparison, it is shown that flux \eqref{llf_u} is entropy stable provided $Q_{i+\frac{1}{2}} - Q*_{i+\frac{1}{2}}\geq 0$, which holds true under condition \eqref{alpha_cond}}. Examples of viscosity coefficient of such entropy stable fluxes are,
\begin{itemize}
	\item $\alpha_{i+\frac{1}{2}}= \displaystyle \max_{\lambda,|\xi|\leq 1/2}\left|\lambda\left[\mathbf{A(u(v}_{i+\frac{1}{2}}(\xi)))\right]\right| $ viz Rusanov viscosity \cite{Rusanov}.
	\item $\alpha_{i+\frac{1}{2}} = \displaystyle \max_{\lambda,\mathbf{u}}|\lambda\left[\mathbf{A(u)}\right]|$ viz Lax-Friedrichs viscosity \cite{Friedrichs}. 
\end{itemize}
\hfill $\square$
\pagebreak
\subsection{High order entropy stable TeCNO flux}\label{sec5b}
The high order entropy stable numerical fluxes in \cite{Fjordholm2012,BbRk,DUAN2021110136,LIU2019104266} etc can be written as
\begin{equation}\label{tecno_u}
	{\mathbf{F}}_{i+\frac{1}{2}}=\mathbf{F}_{i+\frac{1}{2}}^*-\frac{1}{2} \mathbf{\tilde{R}}_{i+\frac{1}{2}}\mathbf{\tilde{\Lambda}}_{i+\frac{1}{2}}\left<\left<\mathbf{\tilde{w}}\right>\right>_{i+\frac{1}{2}}.
\end{equation}
On using scaled entropy variable relation \eqref{tildew}, equation \eqref{tecno_u} can be written as 
\begin{equation}\label{eq_new}
	{\mathbf{F}^*}_{i+\frac{1}{2}}-{\mathbf{F}}_{i+\frac{1}{2}}=\frac{1}{2} \mathbf{\tilde{R}}_{i+\frac{1}{2}}\tilde{\mathbf{\Lambda}}_{i+\frac{1}{2}} \mathbf{\tilde{R}}_{i+\frac{1}{2}}^T\left<\left<\mathbf{v}\right>\right>_{i+\frac{1}{2}},\,
\end{equation}
Let flux \eqref{eq_new} satisfies the component wise flux sign stability condition \eqref{signpropu} i.e., 	\begin{equation}\label{signpropu_here}
	sign\left(\mathbf{F}^*_{i+\frac{1}{2}}-{\mathbf{F}}_{i+\frac{1}{2}}\right)= 	sign\left(\left[\left[\mathbf{v}\right]\right]_{i+\frac{1}{2}}\right)
\end{equation}
which holds provided
$$sign\left(\left<\left<\mathbf{v}\right>\right>_{i+\frac{1}{2}}\right)=sign\left(\left[\left[\mathbf{v}\right]\right]_{i+\frac{1}{2}}\right).$$

Thus by Lemma \ref{lem-new}, $ \left<\left<\mathbf{v}\right>\right>_{i+\frac{1}{2}}= \mathbf{S}\left[\left[\mathbf{v}\right]\right]_{i+\frac{1}{2}}$ and therefore, $$ \left<\left<\mathbf{v}\right>\right>_{i+\frac{1}{2}}= (\mathbf{\tilde{R}}_{i+\frac{1}{2}}^T)^{-1}\mathbf{\tilde{R}}_{i+\frac{1}{2}}^T\mathbf{S}\left[\left[\mathbf{v}\right]\right]_{i+\frac{1}{2}}\;\text{or}\;
 \mathbf{\tilde{R}}_{i+\frac{1}{2}}^T\left<\left<\mathbf{v}\right>\right>_{i+\frac{1}{2}}= \mathbf{S}\left(\mathbf{\tilde{R}}_{i+\frac{1}{2}}^T\left[\left[\mathbf{v}\right]\right]_{i+\frac{1}{2}}\right).$$
Since the scaling matrix $\mathbf{S} = diag(s_{11},s_{22}, \dots, s_{mm}), s_{ii}\geq 0,\, \text{for}\, i=1,2,\dots m$. Therefore flux sign stability condition \eqref{signpropu_here} reduces to following component wise sign stability condition  
\begin{equation} \label{RTv}
	sign\left(\mathbf{\tilde{R}}_{i+\frac{1}{2}}^T\left<\left<\mathbf{v}\right>\right>_{i+\frac{1}{2}}\right)= sign\left(\mathbf{\tilde{R}}_{i+\frac{1}{2}}^T\left[\left[\mathbf{v}\right]\right]_{i+\frac{1}{2}}\right).
\end{equation}
Thus flux sign stability satisfying flux \eqref{tecno_u} is entropy stable provided \eqref{RTv} holds. Note that condition \eqref{RTv} is shown for the entrpy stability of TeCNO flux in \cite{Fjordholm2012}, that.

\hfill $\square$\\
\subsection{Example of a flux sign stable pair $(\mathbf{F}^{*,4}, \mathbf{F}^{s,3})$}
 In order to show that high order flux sign stable pairs exists, this section particularly shows that forth order entropy conservative and and third non-oscillatory WENO flux pair satisfies the flux sign stability property across discontinuities. Note that the entropy changes across discontinuities \cite{Ismail2009} therefore, it is sufficient to achieve entropy stability around discontinuous solution region. Note that the flux sign stability property is defined component wise therefore following proof in scalar setting suffice.  We note that discontinuities in the solution can be characterized by locally significant jump (LSJ) \cite{BbRk}. More precisely, let $\{u_{i}\}$ denotes the set of discrete data values of a function $u(x)$ at $x_{i}$ then the discrete set $\{u_{i}\}$ has a locally significant jump in the interval $[x_{i},x_{i+1}]$ if 
\begin{equation}
|u_{i+1}-u_{i}|>max(|u_{i}-u_{i-1}|,|u_{i+2}-u_{i+1}|). \label{lbj}
\end{equation}
i.e., jump in the interval $[x_{i},x_{i+1}]$ is bigger than the jumps on immediate left and right intervals. Let $r_{i}^+=\frac{[[ u]]_{i-\frac{1}{2}}}{[[ u]]_{i+\frac{1}{2}}}$ and $r_{i}^-=\frac{[[ u]]_{i+\frac{1}{2}}}{[[ u]]_{i-\frac{1}{2}}}$, 
where $[[ * ]]_{i+\frac{1}{2}}=(*)_{i+1}-(*)_{i} $.
It follows from \eqref{lbj} that for a discontinuity present in cell $[x_i, x_{i+1}]$, conditions $|[[ u ]]_{i-1/2}|<|[[ u ]]_{i+1/2}|$ and  $|[[ u ]]_{i+3/2}|<|[[ v ]]_{i+1/2}|$ hold which in turn implies that at discontinuity 
\begin{equation}\label{modR}
|r_{i}^+| <1\; \mbox{and}\;|r_{i+1}^-| < 1.
\end{equation}
\begin{theorem} Let $F^{*,4} $ be the forth order entropy conservative flux  and $F^{s,3}$ be third WENO flux then for the linear problem \begin{equation}\label{linprob}
	u_t + u_x=0, u(x,0) = u_0(x),
	\end{equation} 
	the condition 
	\begin{equation}\label{EC4W3}
	sign(F^{*,4}_{i+\frac{1}{2}}- F^{s,3}_{i+\frac{1}{2}}) = sign([[{v}]]_{i+1/2} )
	\end{equation} is satisfied at solution discontinuities.
\end{theorem} 
{\bf Proof:} 
For scalar problem \eqref{linprob} $v=u$ and 
\begin{equation}\label{ec4}
F^{*,4} = \frac{4}{3}\left(\frac{u_i+u_{i+1}}{2}\right) -\frac{1}{6}\left(\frac{u_{i-1}+u_{i+1}}{2}+ \frac{u_{i}+u_{i+2}}{2}\right)
\end{equation} 
\begin{equation}\label{w3}
F^{s,3} = w_0\left(\frac{3}{2}u_i - \frac{1}{2}u_{i-1}\right) + w_1\left(\frac{1}{2}(u_i + u_{i+1})\right)
\end{equation}
where $w_0, w_1$ are non-linear WENO weights satisfying $w_0+w_1 =1, w_0\geq 0, w_1\geq 0$. Then we have
\begin{eqnarray}
F^{*,4} -F^{s,3} &=& \frac{w_0}{2}\left(u_{i+1} - u_i\right) + \left(\frac{1}{12}-\frac{w_0}{2}\right)(u_i-u_{i-1}) -\frac{1}{12}(u_{i+2}-u_{i+1})\nonumber\\
&=&\left(\frac{w_0}{2}(1-r_{i}^+) + \frac{1}{12}r_{i}^+ -\frac{1}{12}r_{i+1}^-\right)[[ u]]_{i+\frac{1}{2}}\nonumber\\
&=& D[[ u]]_{i+\frac{1}{2}}
\end{eqnarray}	
where $D =\frac{w_0}{2}(1-r_{i}^+) + \frac{1}{12}r_{i}^+ -\frac{1}{12}r_{i+1}^-$. Clearly at solution discontinuity, condition \eqref{EC4W3} will satisfy provided $D\geq 0$ i.e., 
\begin{equation}\label{inq1}
(1-r_{i}^+) \geq  \frac{1}{6w_0}(r_{i+1}^- -r_{i}^+)
\end{equation}
We refer \cite{Sabana}, where third order WENO weights $w_0$ for non-oscillatory approximation is characterized as a function of $r_i^+$. It is shown therein  that $w_0(r_i^+) \geq 1/3$ for $|r_i^+|<1$. Thus the above inequality \eqref{inq1} satisfies at discontinuity if
\begin{equation}\label{inq2}
(1-r_{i}^+) \geq  \frac{1}{2}(r_{i+1}^- -r_{i}^+),
\end{equation}
Or 
\begin{equation}\label{inq2}
(1-r_{i}^+) \geq  (r_{i+1}^- -r_{i}^+) \Rightarrow r_{i+1}^- \leq 1,
\end{equation}
which is indeed the case at LSJ or discontinuities. Therefore, at discontinuities flux sign stability condition \eqref{EC4W3} is satisfied. Thus at discontinuities the entropy stable flux $\hat{F}$ in \eqref{ES_Flux} constructed by \eqref{finalESflux} reduces to non-oscillatory $F^{s,3}$ flux i.e., formulation \eqref{finalESflux} yields a non-oscillatory entropy stable scheme.
 
\section{Numerical Result}\label{sec6}
The approach to construct entropy stable flux \eqref{finalESflux} is generic and robust as it can work with richly available any entropy conservative $\mathbf{F}^*$ and any non-oscillatory $\mathbf{F}^s$ fluxes. Entropy conservative fluxes  for large class of hyperbolic systems can be found in  \cite{LeFloch2002,Ismail2009,Fjordholm2012,DUAN2021110136,Chandrashekar2013,Winters201672} whereas non-oscillatory fluxes can be chosen from  \cite{Shu1997,liu1994weighted,biswas2018accuracy,jiang1996efficient,Sabana,HARTEN1983357,ZhangTVD,Zhang2011,borges2008improved}. In the following, several numerical tests are presented to demonstrate the accuracy and non-oscillatory capability of entropy stable schemes using various combinations. Precisely, we choose following setting to construct entropy stable non-oscillatory fluxes and numerical computations: \\
  {\bf Entropy conservative fluxes:}  We choose entropy conservative flux $\mathbf{F}^*$ given in \cite{Tadmor1987, LeFloch2002}.
  In particular, we use following second, fourth and sixth order entropy conservative fluxes
  \begin{equation}\label{ECflux4}
  \mathbf{F}^{*,4}_{i+\frac{1}{2}} = \frac{4}{3}\mathbf{F}_{i+\frac{1}{2}}^*(\mathbf{u}_i,\mathbf{u}_{i+1}) -\frac{1}{6}(\mathbf{F}_{i+\frac{1}{2}}^*(\mathbf{u}_{i-1},\mathbf{u}_{i+1}) + \mathbf{F}_{i+\frac{1}{2}}^*(\mathbf{u}_{i},\mathbf{u}_{i+2}))
  \end{equation}

  \begin{eqnarray}\label{ECflux6}
  \mathbf{F}^{*,6}_{i+\frac{1}{2}} = \frac{3}{2}\mathbf{F}_{i+\frac{1}{2}}^*(\mathbf{u}_i,\mathbf{u}_{i+1}) -\frac{3}{10}(\mathbf{F}_{i+\frac{1}{2}}^*(\mathbf{u}_{i-1},\mathbf{u}_{i+1}) + \mathbf{F}_{i+\frac{1}{2}}^*(\mathbf{u}_{i},\mathbf{u}_{i+2})) \nonumber \\ 
  +
  \frac{1}{30}(\mathbf{F}_{i+\frac{1}{2}}^*(\mathbf{u}_{i-2},\mathbf{u}_{i+1}) + \mathbf{F}_{i+\frac{1}{2}}^*(\mathbf{u}_{i-1},\mathbf{u}_{i+2}) + \mathbf{F}_{i+\frac{1}{2}}^*(\mathbf{u}_{i},\mathbf{u}_{i+3}))
  \end{eqnarray}
  where $\displaystyle \mathbf{F}_{i+\frac{1}{2}}^*(\mathbf{a}, \mathbf{b})$ is second order entropy conservative fluxes satisfying \eqref{eq_ec}  \cite{Fjordholm2012,DUAN2021110136,Ismail2009,Tadmor2003,Tadmor2016}.\\
 {\bf Non-oscillatory fluxes:} The used non-oscillatory fluxes $\mathbf{F}^s$ are first order Local Lax-Friedrichs or Rusanov flux \eqref{llf_u}, second and third order ENO fluxes \cite{Shu1997}, third and fifth order WENO-JS fluxes with Jiang-Shu weights  \cite{jiang1996efficient,Sabana}. For comparison of results in Figure \ref{LintestIC2c} WENO-Z fluxes with weight \cite{borges2008improved} is also used. In the $m$-system case \eqref{conservation_laws}, a straightforward and computationally efficient component-wise reconstruction \cite{Shu2009} is used as follows:
  \begin{enumerate}
  	\item Split the physical flux  $\mathbf{f(u)}$ in \eqref{conservation_laws} as, 
  	\begin{equation}
  	\mathbf{f(u)}=\mathbf{f(u)}^++\mathbf{f(u)}^-,
  	\end{equation}
  	such that $\dfrac{d\mathbf{f(u)}}{d\mathbf{u}}^+\geq 0$ and $\dfrac{d\mathbf{f(u)}}{d\mathbf{u}}^-\leq 0$. 
  	\item The $m^{th}$ component of discrete numerical flux ${{\mathbf{F}}}_{i+\frac{1}{2}}^{\pm}$ using ENO/WENO reconstructions is defined as, \begin{equation}
  	{{F}^{(m)\pm}}_{i+\frac{1}{2}}=v_{i+\frac{1}{2}}^{\mp},
  	\end{equation}
  	where $\bar{v}_{i}=f^{(m)\pm}_i$ is taken as input.
  	\item Then the non-oscillatory ENO/WENO  numerical flux is obtained by, 
  	\begin{equation}
  	{{\mathbf{F}}}_{i+\frac{1}{2}}={{\mathbf{F}}}_{i+\frac{1}{2}}^{+}+{{\mathbf{F}}}_{i+\frac{1}{2}}^{-}.
  	\end{equation}
  \end{enumerate}
 {\bf Time integration:} For time marching in \eqref{semi_scheme}, following explicit third order strong stability preserving(SSP) Runge-Kutta methods is used \cite{gottlieb1998total,shu1988total}.   
 \begin{subequations} \label{tvdRK3}
 	\begin{equation}
  \mathbf{u}^{(1)}=\mathbf{u}^n+\triangle t \mathcal{L}(\mathbf{u}^n),
  \end{equation}
\begin{equation}  \mathbf{u}^{(2)}=\frac{3}{4} \mathbf{u}^{(n)}+\frac{1}{4}\left(\mathbf{u}^{(1)}+\triangle t \mathcal{L}(\mathbf{u}^{(1)}\right), 
	\end{equation}
\begin{equation} 
  \mathbf{u}^{(n+1)}=\frac{1}{3}+\frac{2}{3} \mathbf{u}^{(2)}+\frac{2}{3}\triangle t \mathcal{L}(\mathbf{u}^{(2)}),
  \end{equation}
  \end{subequations} 
  where
  $ [\mathcal{L}(\mathbf{u})]_i = -\frac{1}{\Delta x}\left[\left[\mathbf{\hat{F}}\right]\right]_{i}$. The extension of the schemes for multi-dimensional test problems on regular mesh is done by dimension by dimension or tensorial approach  \cite{Shu1997, Fjordholm2012} is utilized. Following the approach in \cite{rathan2018modified}, only for accuracy tests, time step is taken as $\Delta t = (\Delta x)^5/3$ so that the third order method \eqref{tvdRK3} in time is effectively fifth-order. In all other test cases time step is taken such that it satisfies $ \frac{\max_{u}|f'(u)| \Delta t }{\Delta x}=CFL< 1.$   
{\bf Name Convention:} The following name convention is used to annotate the results: legend {\it EC-m-$F^w$-n} denotes the numerical result using entropy stable scheme constructed by combination of $m^{th}$ order accurate entropy conservative flux $F^*$ and $n^{th}$ order non-oscillatory flux $F^s$. For example, EC-6-WENOJS-5 represents entropy stable flux obtained by combining $6^{th}$ order entropy conservative flux \eqref{ECflux6} with the fifth order non-oscillatory WENOJS flux \cite{borges2008improved}.
 
\subsection{Scalar conservation law} We consider scalar problems to analyze the accuracy and non-oscillatory nature of various entropy stable schemes obtained by different combinations of entropy conservative fluxes $F^*$ and non-oscillatory fluxes $F^s$.
\subsubsection{\bf Linear transport equation:} 
\begin{figure}
	\begin{tabular}{cc}
		\includegraphics[scale =0.55]{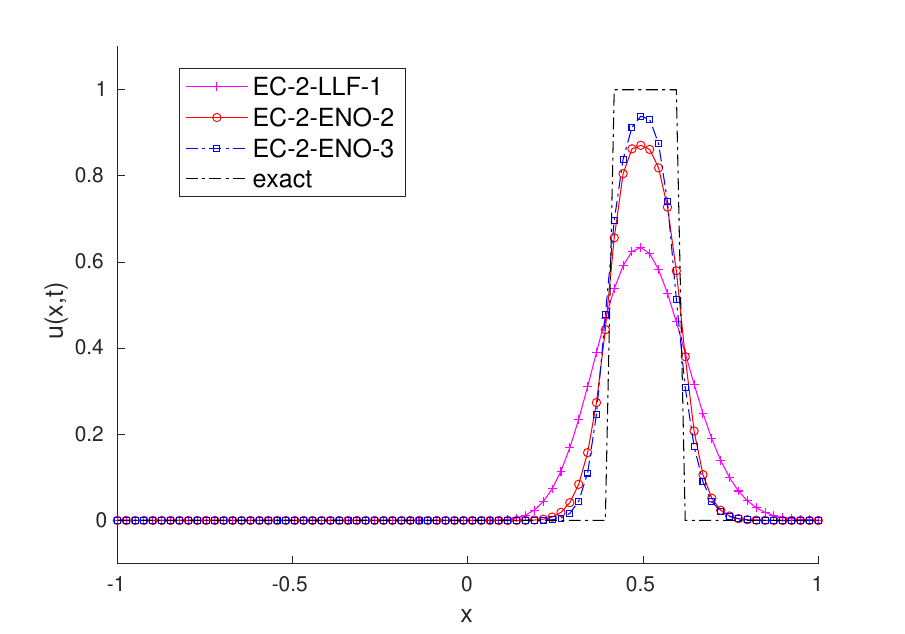} &
		\includegraphics[scale =0.55]{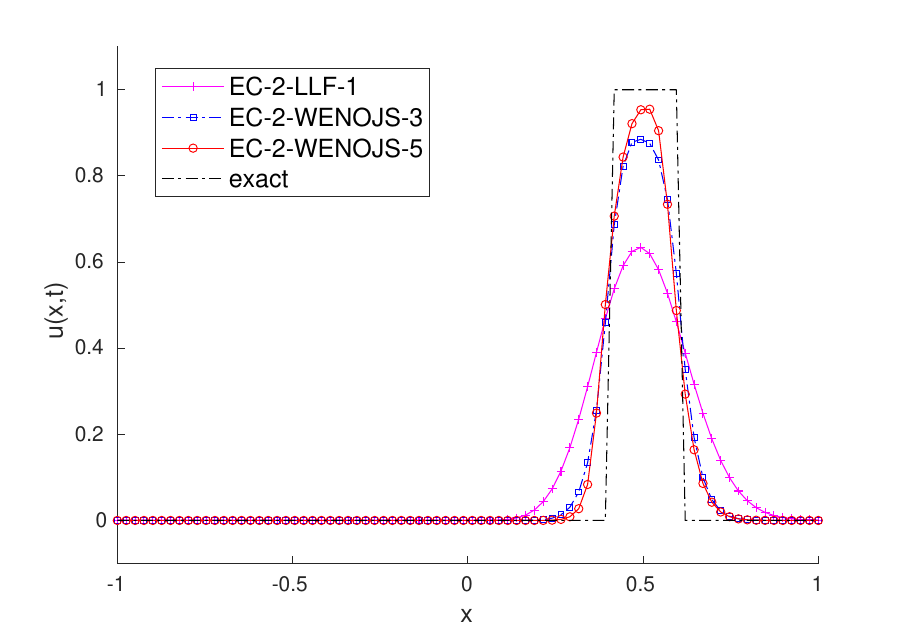}\\
			\includegraphics[scale =0.55]{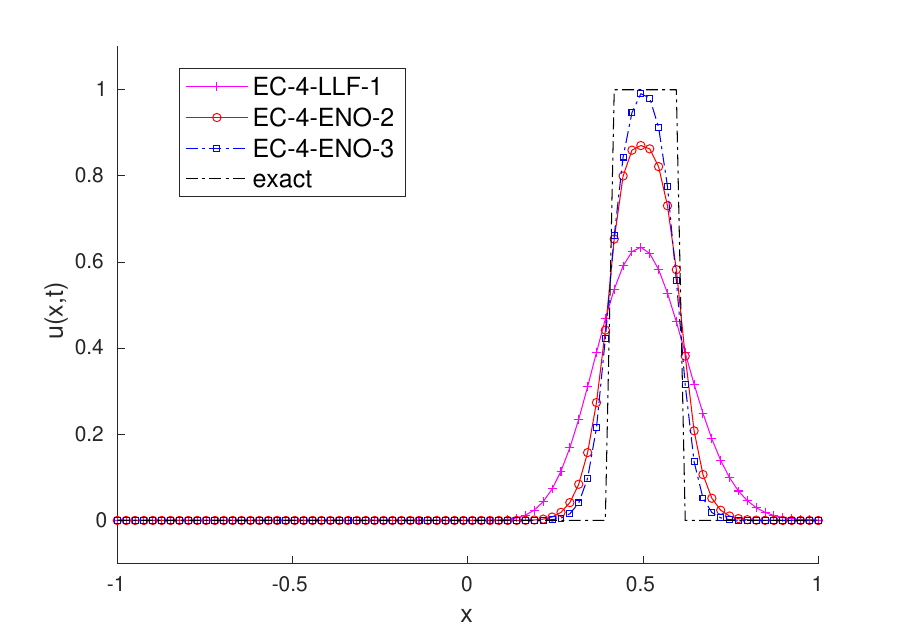} &
		\includegraphics[scale =0.55]{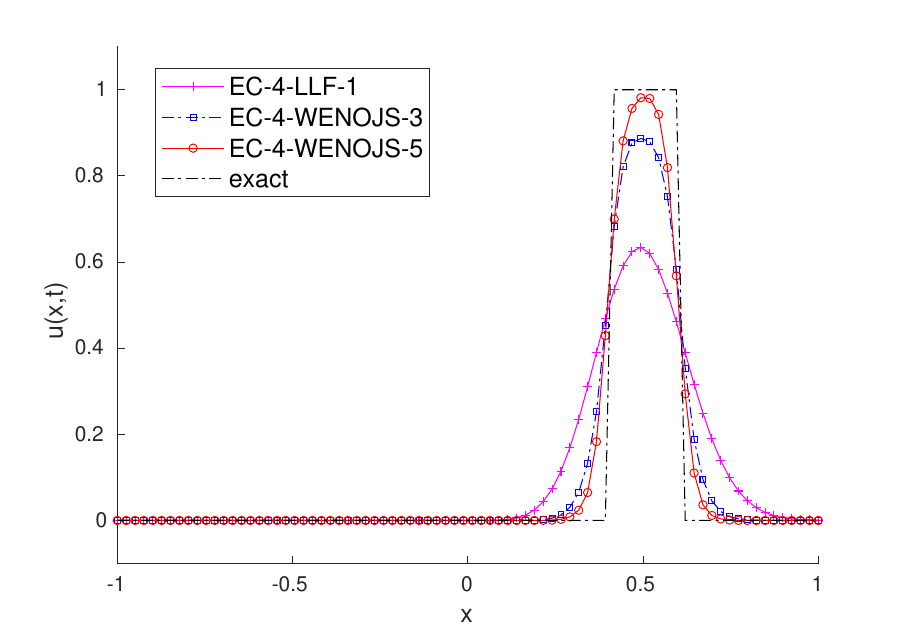}\\
			\includegraphics[scale =0.55]{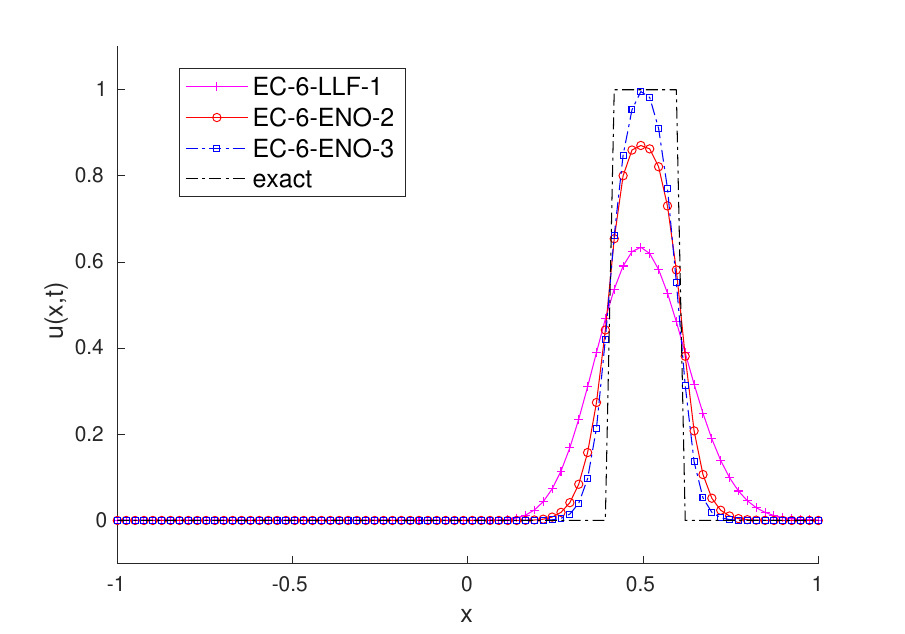} &
		\includegraphics[scale =0.55]{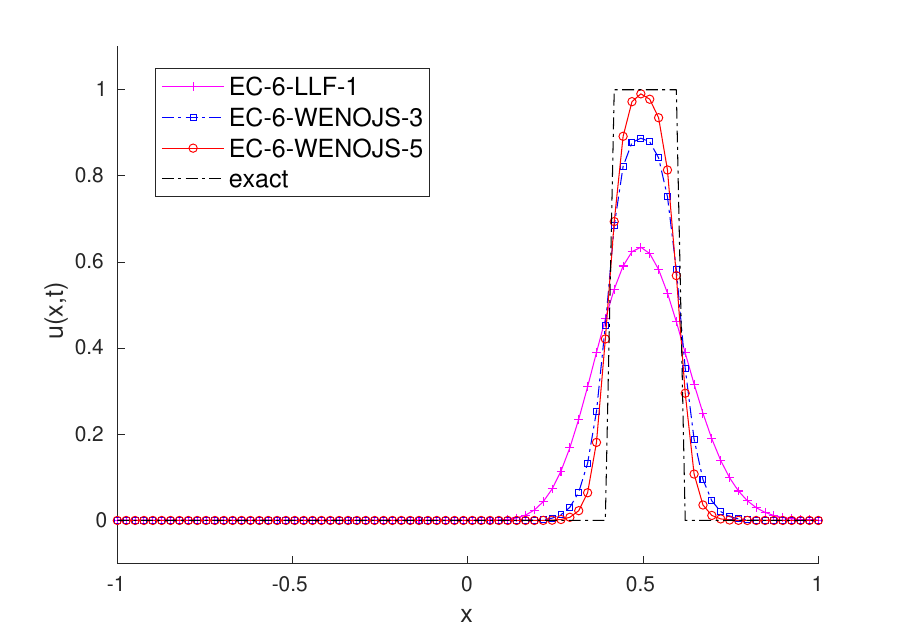}	
	\end{tabular}
\caption{\label{LintestIC2a} {\it Linear transport equation for non-oscillatory property:}  Each sub figure corresponds to the fixed EC flux $F^*$ with different $F^s$ fluxes.}
\end{figure}

\begin{figure}
	\begin{tabular}{cc}
		\includegraphics[scale =0.55]{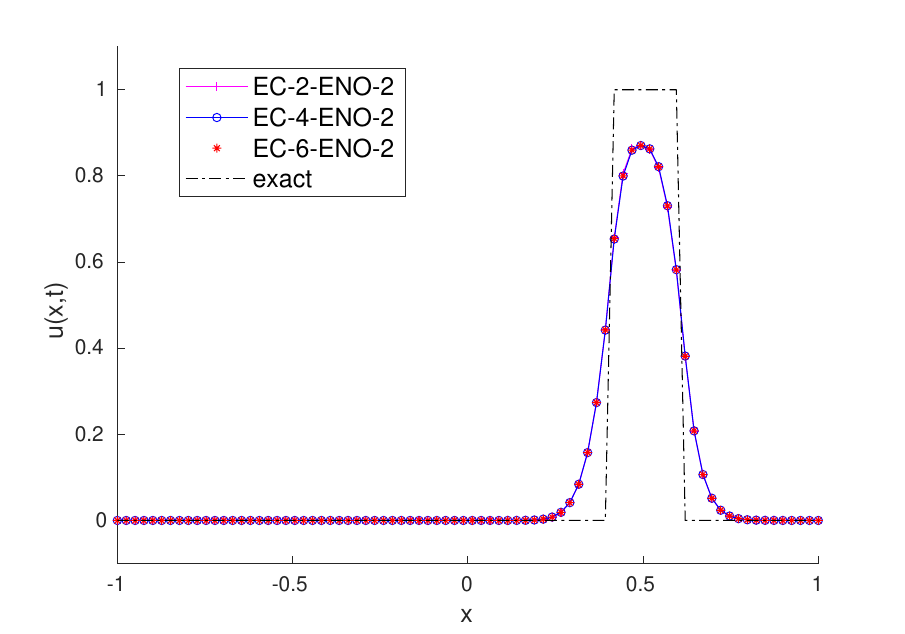}
		 & \includegraphics[scale =0.55]{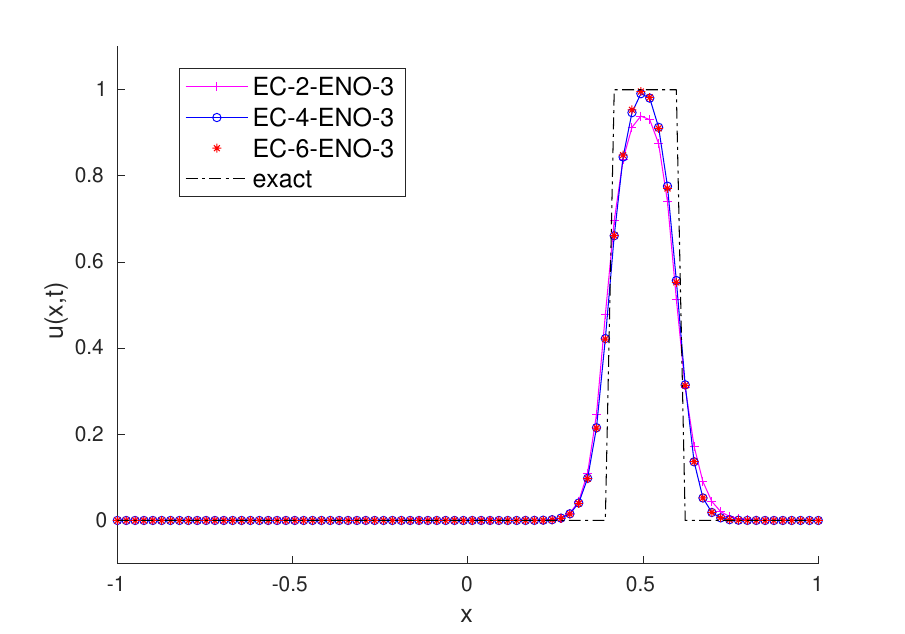}\\
	\includegraphics[scale =0.55]{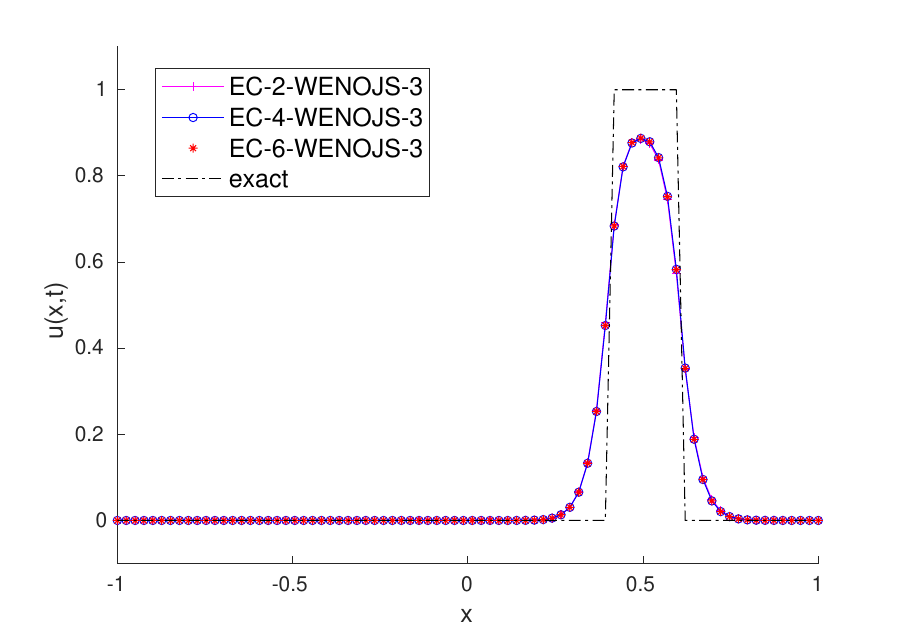} & 
		 \includegraphics[scale =0.55]{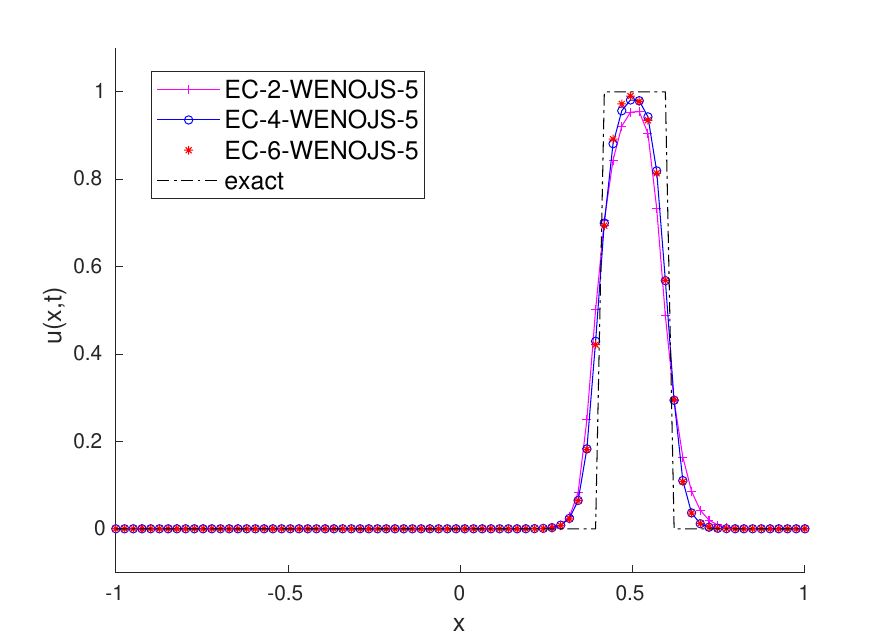}\\
	\end{tabular}
	\caption{\label{LintestIC2b} {\it Linear transport equation for non-oscillatory property:} Each sub figure corresponds to different EC flux with fixed $F^s$ flux.}
\end{figure}

\begin{figure}
	\begin{tabular}{cc}
	\includegraphics[scale=0.55]{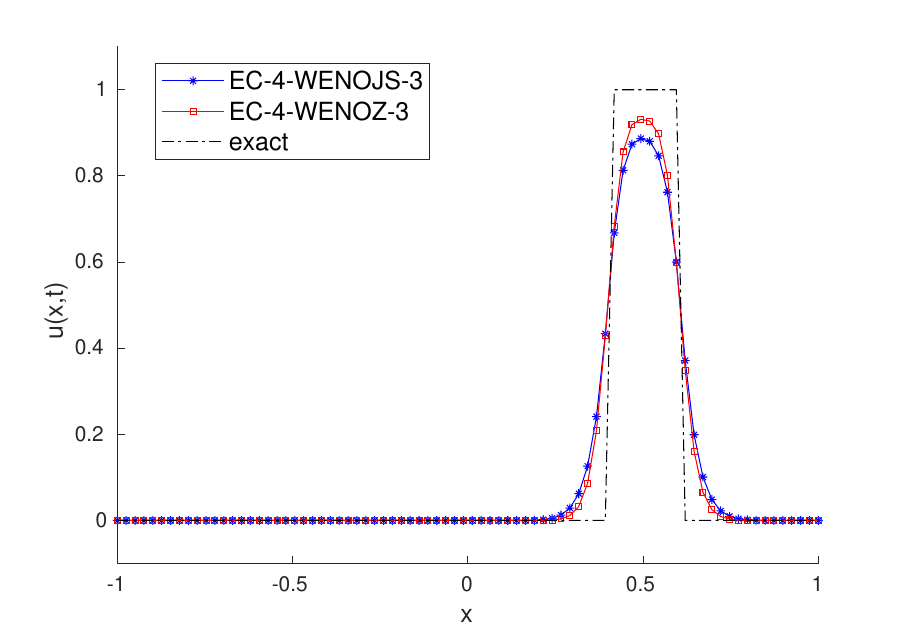}
	&	\includegraphics[scale=0.55]{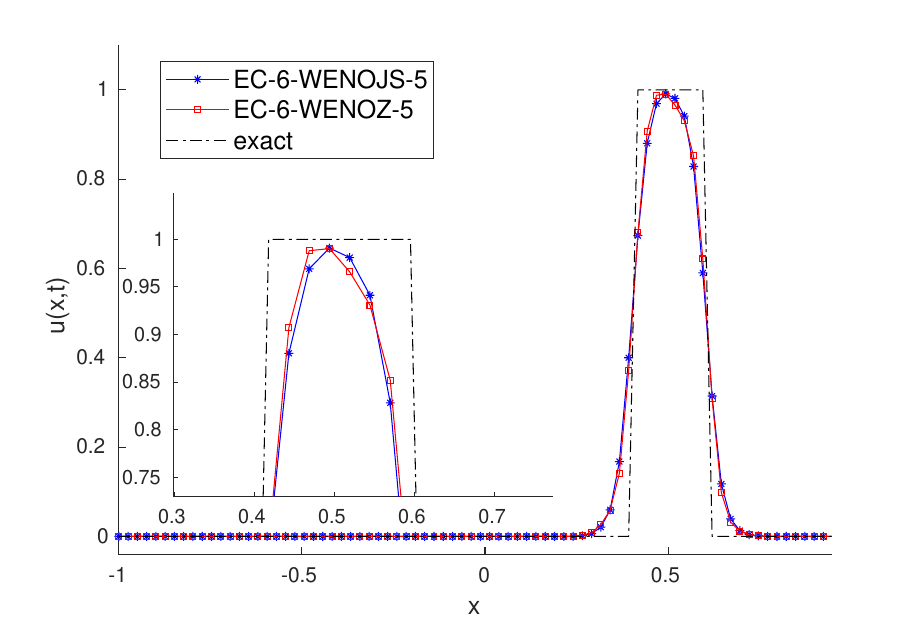}\\
	\end{tabular}
	\caption{\label{LintestIC2c} {\it Linear transport equation for non-oscillatory property:} Different EC fluxes with WENOJS and WENOZ fluxes.}
	\label{fig:lintest4ec-4weno-compare}
\end{figure}

Consider the following linear equation 
  \begin{equation}\label{transporteq}
  u_t+u_x=0, 
  \end{equation}
  along with the following initial conditions.
  \begin{itemize}
  	\item{IC1: To test accuracy of schemes:}  \begin{equation}\label{LinearIC1}
  	u(x,0)=-sin(\pi x) ,\;\;\; x\in [-1,1].  
  	\end{equation}
  	In Table \ref{tab1LinearIC1} to \ref{tab4LinearIC1}, the $L^\infty$ and $L^1$ convergence rate of different non-oscillatory entropy stable schemes is given and compared with the corresponding non-oscillatory scheme with respect to initial condition \eqref{LinearIC1}. It can be seen that the entropy stable scheme using $m^{th}$ order entropy conservative flux and $n^{th}$ order non-oscillatory flux {\it EC-m-$F^w$-n} maintains formal order of accuracy i.e., show the $k^{th}$ order sub-optimal convergence rate, where $k=\min(m,n)$. .
  	\item{IC2: To test non-oscillatory property of schemes:} \begin{equation}\label{LinearIC2}
  	u(x,0)=\left\{\begin{array}{cc}
  	1 & |x|\leq \frac{1}{10}\\
  	0 & else
  	\end{array}\right.,\;\;\;x\in [-1,1]. \\
  	\end{equation}  
  	In figures \ref{LintestIC2a}, \ref{LintestIC2b} and \ref{LintestIC2c}, computed solutions at time $T_f=0.5$ corresponding to initial condition \eqref{LinearIC2} by various entropy stable schemes are given and compared. It can be observed that these entropy stable schemes do not produce any \textcolor{blue}{nonphysical spurious oscillations} in the vicinity of discontinuities. From Figures \ref{LintestIC2a} and \ref{LintestIC2b} it can be observed that the resolution of discontinuities by entropy stable schemes is characterized by underlying $k^{th}$ order flux $F^*=EC${\it-m} or $F^s=F^w${\it-n}. 
  	In figure \ref{LintestIC2c}, results are compared between entropy stable schemes using WENO-JS weights \cite{jiang1996efficient} and their modified version improved WENO-Z weights \cite{borges2008improved} to achieve optimal third-order accuracy regardless of the critical point. This improvement between third order WENO-JS and WENO-Z fluxes also reflects clearly in the computational results in Figure \ref{LintestIC2c} (Left).   	
  \end{itemize}   	 
\begin{table}[htb!]
	\centering	
	\begin{tabular}{|cc|}
		\hline
		{\bf ENO-2} & {\bf EC2-ENO-2}\\
		\begin{tabular}{|c|c|c|c|c|}
			\hline N  & $L^\infty$ error  &  Rate &    $L^1$ error    & Rate \\
			\hline 20  & 1.0365e-01 &        - &  1.0093e-01 &        -\\
			\hline 40  & 4.1705e-02 &  1.31 &  2.8773e-02 &  1.81\\
			\hline 80 &  1.7092e-02 &  1.29 &  7.8886e-03  & 1.87\\
			\hline 160 &  6.9733e-03 &  1.29 & 2.0813e-03  & 1.92\\
			\hline 320 &  2.8250e-03 &  1.30 &  5.4724e-04  & 1.93\\
			\hline 640 &  1.1408e-03 &  1.31 & 1.4310e-04 &  1.94\\
			\hline 
		\end{tabular}&		
		\begin{tabular}{|c|c|c|c|c|}
			\hline N  & $L^\infty$ error  &  Rate &    $L^1$ error    & Rate \\
			\hline 20  & 1.0365e-01   &      - &  1.0093e-01 &   -\\
			\hline 40  & 4.1705e-02  & 1.31  & 2.8773e-02  & 1.81\\
			\hline 80  & 1.7092e-02  & 1.29  & 7.8886e-03  & 1.87\\
			\hline 160 & 6.9733e-03  & 1.29  & 2.0813e-03  & 1.92\\
			\hline 320 & 2.8250e-03  & 1.30  & 5.4724e-04  & 1.93\\
			\hline 640 & 1.1408e-03  & 1.31  & 1.4310e-04 &  1.94\\
			\hline				
		\end{tabular}\\ & \\
		{\bf EC4-ENO-2} & {\bf EC6-ENO-2}\\
		\begin{tabular}{|c|c|c|c|c|}
			\hline N  & $L^\infty$ error  &  Rate &    $L^1$ error    & Rate \\
			\hline 20 & 1.0115e-01 & -  & 1.0057e-01 &  -\\
			\hline 40 & 4.1673e-02  & 1.28  & 2.8451e-02  & 1.82\\
			\hline 80 & 1.7076e-02  & 1.29  & 7.8348e-03  & 1.86\\
			\hline 160 & 6.9579e-03  & 1.30  & 2.0692e-03  & 1.92\\
			\hline 320 & 2.8226e-03  & 1.30  & 5.4484e-04  & 1.93\\
			\hline 640 & 1.1403e-03  & 1.31  & 1.4266e-04  & 1.93\\ 
			\hline				
		\end{tabular}&
		\begin{tabular}{|c|c|c|c|c|}
			\hline N  & $L^\infty$ error  &  Rate &    $L^1$ error    & Rate \\
			\hline 20  & 1.0072e-01  &-1.000 &  1.0025e-01  &-\\
			\hline	40  & 4.1641e-02 &  1.27 &  2.8444e-02 &  1.82\\
			\hline80  & 1.7050e-02  & 1.29  & 7.8323e-03  & 1.86\\
			\hline 160 &  6.9540e-03  & 1.29  & 2.0695e-03 &  1.92\\
			\hline 320 &  2.8217e-03  & 1.30  & 5.4487e-04 &  1.93\\
			\hline640 &  1.1399e-03  & 1.31  & 1.4267e-04 &  1.93\\
 			\hline
		\end{tabular}\\ & \\			
		\hline
	\end{tabular}
	\caption{\label{tab1LinearIC1} {\it Linear transport equation:} Convergence rate of base non-oscillatory schemes and corresponding entropy stable schemes {\it EC-m-$F^s$-n} for initial condition \eqref{LinearIC1}, $\Delta t = \Delta x^{\frac{5}{3}},\; T_f =0.5$.}
\end{table}

\begin{table}[htb!]
	\centering	
	\begin{tabular}{|cc|}
		\hline
		{\bf ENO-3} & {\bf EC2-ENO-3}\\
		\begin{tabular}{|c|c|c|c|c|}
			\hline N  & $L^\infty$ error  &  Rate &    $L^1$ error    & Rate \\
			\hline 20   &   5.4065e-03   &  -   &   6.7349e-03   &  -\\
			\hline 40   &   6.1630e-04   &   3.13   &   7.3093e-04   &  3.20\\
			\hline 80   &   7.2521e-05   &   3.09   &   8.5989e-05   &  3.09\\
			\hline 160   &   8.7912e-06   &   3.04   &   1.0394e-05   &  3.05\\
			\hline 320   &   1.0772e-06   &   3.03   &   1.2804e-06   &  3.02\\
			\hline 640   &   1.3245e-07   &   3.02   &   1.5887e-07   &  3.01\\
			\hline 
		\end{tabular}&		
		\begin{tabular}{|c|c|c|c|c|}
			\hline N  & $L^\infty$ error  &  Rate &    $L^1$ error    & Rate \\
			\hline   20   &   5.0543e-02   &  -   &   4.7520e-02   &  -   \\
			\hline   40  &   2.1967e-02   &   1.20   &   1.3347e-02   &   1.83   \\
			\hline   80   &   9.4715e-03   &   1.21   &   3.6014e-03   &   1.89   \\
			\hline   160   &   3.9848e-03   &   1.25   &   9.7373e-04   &   1.89   \\
			\hline   320   &   1.6540e-03   &   1.27   &   2.5970e-04   &   1.91   \\
			\hline   640   &   6.7459e-04   &   1.29   &   6.8309e-05   &   1.93   \\
			\hline				
		\end{tabular}\\ & \\
			{\bf EC4-ENO-3} & {\bf EC6-ENO-3}\\
		\begin{tabular}{|c|c|c|c|c|}
			\hline N  & $L^\infty$ error  &  Rate &    $L^1$ error    & Rate \\
			 \hline   20   &   5.5474e-03   &  -   &   6.7533e-03   &  -   \\
			\hline   40   &   6.3152e-04   &   3.13   &   7.3158e-04   &   3.21   \\
			\hline   80  &   7.4491e-05   &   3.08   &   8.6022e-05   &   3.09   \\
			\hline   160   &   9.0451e-06   &   3.04   &   1.0396e-05   &   3.05   \\
			\hline   320   &   1.1136e-06   &   3.02   &   1.2806e-06   &   3.02   \\
			\hline   640   &   1.3725e-07   &   3.02   &   1.5889e-07   &   3.01   \\
			\hline				
		\end{tabular}&
		\begin{tabular}{|c|c|c|c|c|}
			\hline N  & $L^\infty$ error  &  Rate &    $L^1$ error    & Rate \\
			\hline   20   &   7.0393e-03   &  -   &   8.2049e-03   &  -   \\
			\hline   40   &   8.6398e-04   &   3.03   &   9.7036e-04   &   3.08   \\
			\hline   80   &   9.4006e-05   &   3.20   &   1.0891e-04   &   3.16   \\
			\hline   160   &   1.1206e-05   &   3.07   &   1.3135e-05   &   3.05   \\
			\hline   320   &   1.3546e-06   &   3.05   &   1.6129e-06   &   3.03   \\
			\hline   640   &   1.6523e-07   &   3.04   &   1.9983e-07   &   3.01   \\
			\hline
		\end{tabular} \\ &\\			
		\hline
	\end{tabular}
	\caption{\label{tab2LinearIC1} {\it Linear transport equation case:} Convergence rate of non-oscillatory schemes and corresponding entropy stable schemes {\it EC-m-$F^s$-n}, for initial condition \eqref{LinearIC1}, $\Delta t = \Delta x^{\frac{5}{3}},\; T_f =0.5$. }
\end{table}

\begin{table}[htb!]
	\centering	
	\begin{tabular}{|cc|}
		\hline
		{\bf WENOJS-3} & {\bf EC2-WENOJS-3}\\
		\begin{tabular}{|c|c|c|c|c|}
			\hline N  & $L^\infty$ error  &  Rate &    $L^1$ error    & Rate \\
		\hline   20   &   9.5296e-02   &  -   &   9.7590e-02   &  -   \\
		\hline   40   &   3.7882e-02   &   1.33   &   2.4549e-02   &   1.99   \\
		\hline   80   &   1.4499e-02   &   1.39   &   5.9299e-03   &   2.05   \\
		\hline   160   &   4.4866e-03   &   1.69   &   1.1667e-03   &   2.35   \\
		\hline   320   &   7.9278e-04   &   2.50   &   1.3884e-04   &   3.07   \\
		\hline   640   &   5.0792e-05   &   3.96   &   1.3117e-05   &   3.40   \\
		
			\hline 
		\end{tabular}&		
		\begin{tabular}{|c|c|c|c|c|}
			\hline N  & $L^\infty$ error  &  Rate &    $L^1$ error    & Rate \\
		\hline   20   &   9.6620e-02   &  -   &   9.6364e-02   &  -   \\
		\hline   40   &   3.8751e-02   &   1.32   &   2.5351e-02   &   1.93   \\
		\hline   80   &   1.4978e-02   &   1.37   &   6.2663e-03   &   2.02   \\
		\hline   160   &   4.8798e-03   &   1.62   &   1.3280e-03   &   2.24   \\
		\hline   320   &   1.2733e-03   &   1.94   &   2.2698e-04   &   2.55   \\
		\hline   640   &   3.6191e-04   &   1.81   &   4.2453e-05   &   2.42   \\
		
			\hline				
		\end{tabular}\\ & \\
		{\bf EC4-WENOJS-3} & {\bf EC6-WENOJS-3}\\
		\begin{tabular}{|c|c|c|c|c|}
			\hline N  & $L^\infty$ error  &  Rate &    $L^1$ error    & Rate \\
			\hline   20   &   9.3950e-02   &  -   &   9.6604e-02   &  -   \\
			\hline   40   &   3.7700e-02   &   1.32   &   2.4560e-02   &   1.98   \\
			\hline   80   &   1.4451e-02   &   1.38   &   5.8577e-03   &   2.07   \\
			\hline   160   &   4.4611e-03   &   1.70   &   1.1179e-03   &   2.39   \\
			\hline   320   &   8.2318e-04   &   2.44   &   1.2255e-04   &   3.19   \\
			\hline   640   &   1.7013e-04   &   2.27   &   1.2785e-05   &   3.26   \\
			\hline				
		\end{tabular}&
		\begin{tabular}{|c|c|c|c|c|}
			\hline N  & $L^\infty$ error  &  Rate &    $L^1$ error    & Rate \\
			\hline   20   &   6.0920e-03   &  -   &   6.7959e-03   &  -   \\
			\hline   40   &   6.9358e-04   &   3.13   &   7.3364e-04   &   3.21   \\
			\hline   80   &   8.0943e-05   &   3.10   &   8.6142e-05   &   3.09   \\
			\hline   160   &   9.6558e-06   &   3.07   &   1.0401e-05   &   3.05   \\
			\hline   320   &   1.1683e-06   &   3.05   &   1.2808e-06   &   3.02   \\
			\hline   640   &   1.4158e-07   &   3.04   &   1.5890e-07   &   3.01   \\
			\hline
		\end{tabular} \\ &\\			
		\hline
	\end{tabular}

\caption{\label{tab3LinearIC1} {\it Linear transport equation:} Convergence rate of non-oscillatory schemes and corresponding entropy stable schemes {\it EC-m-$F^s$-n}, for initial condition \eqref{LinearIC1}, $CFL=0.8,\; \Delta t = \Delta x^{\frac{5}{3}},\; T_f =0.5$ }
\end{table}     	
\begin{table}[htb!]
		\centering	
		\begin{tabular}{|c|}
			\hline	
		{\bf WENOJS-5}\\
	\begin{tabular}{|c|c|c|c|c|}
			\hline N  & $L^\infty$ error  &  Rate &    $L^1$ error    & Rate \\
			\hline   20   &   1.0633e-03   &  -   &   1.0675e-03   &  -   \\
			\hline   40   &   2.8858e-05   &   5.20   &   2.7413e-05   &   5.28   \\
			\hline   80   &   8.3877e-07   &   5.10   &   7.7734e-07   &   5.14   \\
			\hline   160   &   2.4196e-08   &   5.12   &   2.3317e-08   &   5.06   \\
			\hline   320   &   6.6638e-10   &   5.18   &   7.1625e-10   &   5.02   \\
			\hline   640   &   1.9484e-11   &   5.10   &   2.2600e-11   &   4.99   \\
			\hline 
		\end{tabular}\\	
	EC4-WENOJS-5\\
	\begin{tabular}{|c|c|c|c|c|}
		\hline   20   &   9.0285e-04   &  -   &   1.0550e-03   &  -   \\
		\hline   40   &   1.0566e-04   &   3.10   &   5.6640e-05   &   4.22   \\
		\hline   80   &   1.1828e-05   &   3.16   &   3.5163e-06   &   4.01   \\
		\hline   160   &   1.3147e-06   &   3.17   &   2.2882e-07   &   3.94   \\
		\hline   320   &   1.4447e-07   &   3.19   &   1.4882e-08   &   3.94   \\
		\hline   640   &   1.5852e-08   &   3.19   &   9.6427e-10   &   3.95   \\
		\hline				
	\end{tabular}\\
		EC6-WENOJS-5\\
		\begin{tabular}{|c|c|c|c|c|}
			\hline   20   &   9.8150e-04   &  -   &   1.0138e-03   &  -   \\
			\hline   40   &   2.7221e-05   &   5.17   &   2.6652e-05   &   5.25   \\
			\hline   80   &   8.0889e-07   &   5.07   &   7.7417e-07   &   5.11   \\
			\hline   160   &   2.3875e-08   &   5.08   &   2.3322e-08   &   5.05   \\
			\hline   320   &   6.6638e-10   &   5.16   &   7.1625e-10   &   5.03   \\
			\hline   640   &   1.9488e-11   &   5.10   &   2.2600e-11   &   4.99   \\ \hline
			\end{tabular}\\ \hline
	\end{tabular}
	\caption{\label{tab4LinearIC1} {\it Linear transport equation:} Convergence rate of non-oscillatory schemes and corresponding entropy stable schemes {\it EC-m-$F^s$-n}, $m>n$, for initial condition \eqref{LinearIC1}, $CFL=0.8,\; \Delta t = \Delta x^{\frac{5}{3}},\; T_f =0.5$ }
\end{table}

\subsubsection{Burgers equation}
\begin{figure}
	\begin{tabular}{cc}
	\includegraphics[scale=0.55]{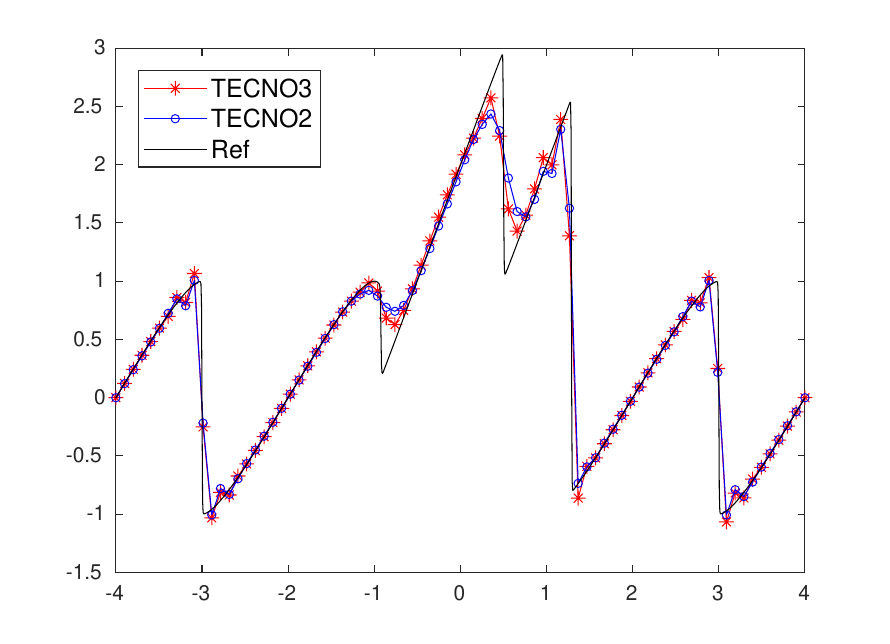} &\includegraphics[scale=0.55]{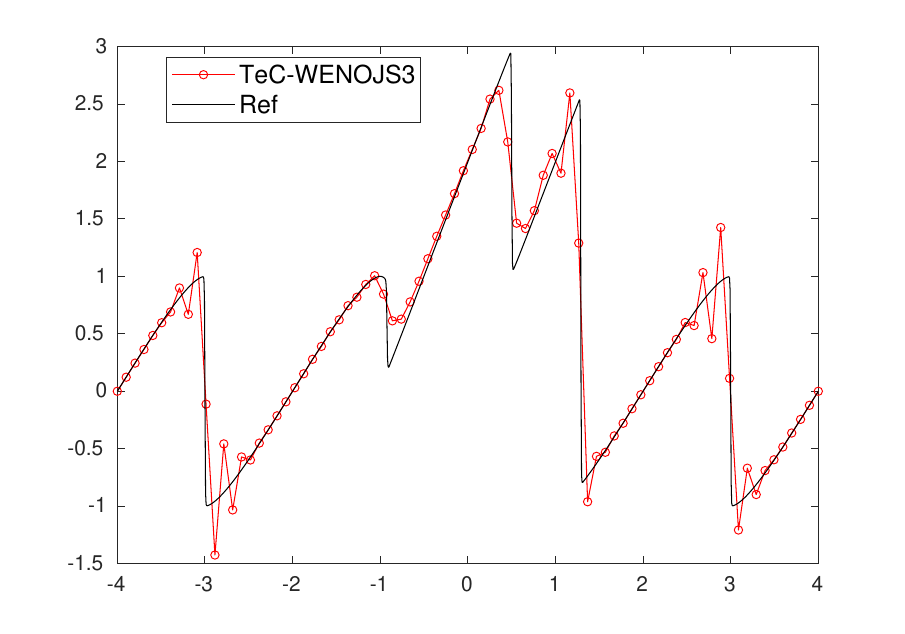}\\
	(a) & (b)\\
	\includegraphics[scale=0.55]{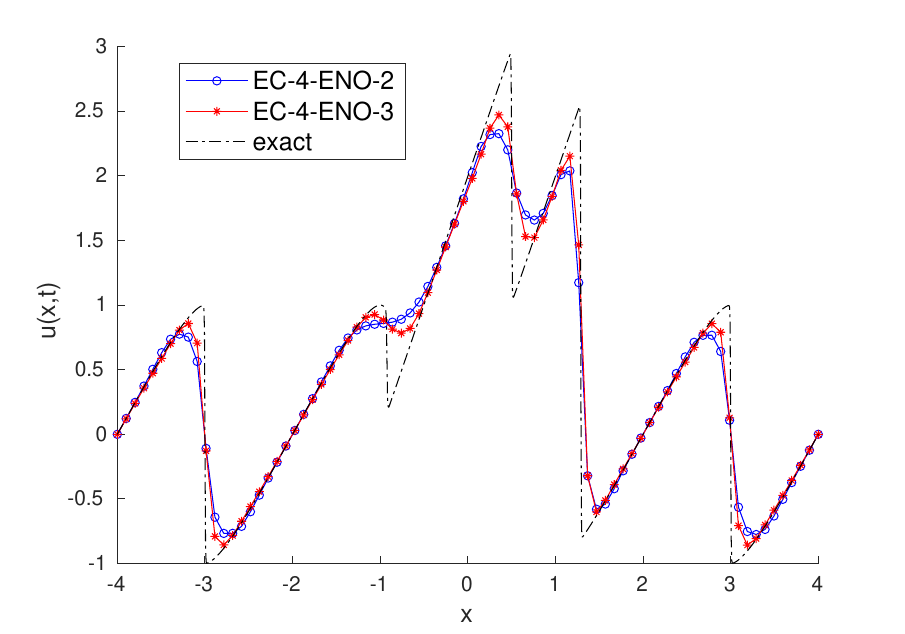} &\includegraphics[scale=0.55]{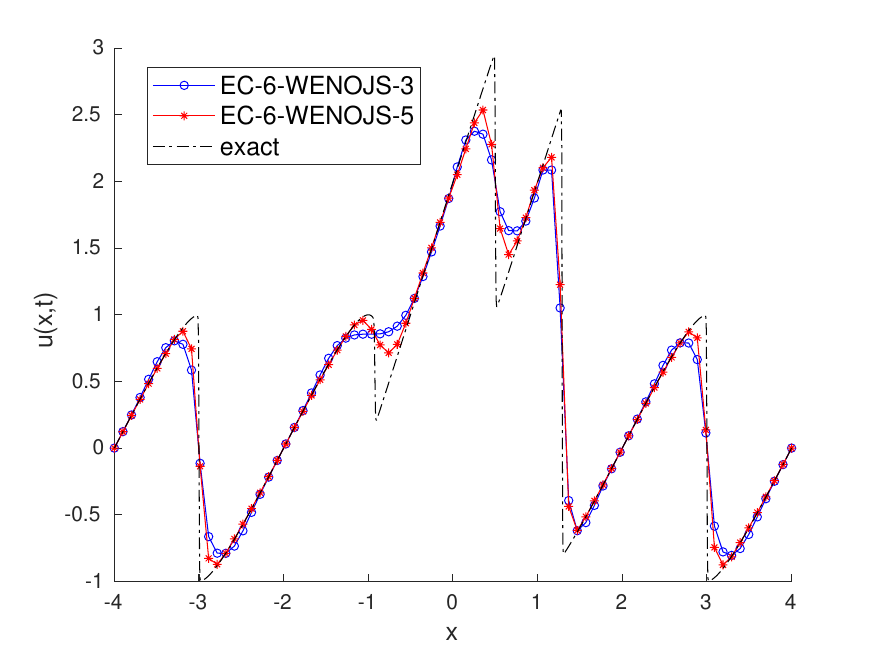}\\
	(c) & (d) 
	\end{tabular}
	\caption{Non-oscillatory solution of Burgers equation by {\it EC-m-$F^w$-n} schemes corresponding to \eqref{BurgerIC2} at $T_f=0.5,\; CFL = 0.8,\; N=80$. Small spurious oscillations can be observed in solution by the TECNO \cite{Fjordholm2012} and TeC-WENOJS3 \cite{BbRk} schemes.}
	\label{fig:burgertest6eceno}
\end{figure}

Consider the Burger's equation \\ \\ \\
\begin{equation} \label{Burgerseq}
	u_t+\bigg(\frac{u^2}{2}\bigg)_x=0,\; x\in [-1,1],
\end{equation}
with periodic boundary condition. The following initial conditions are chosen
\begin{itemize}
	\item{IC3: To test accuracy of schemes:}
\begin{eqnarray}\label{BurgerIC1}
	u(x,0)=1+ \frac{1}{2}sin(\pi x).
	\end{eqnarray}
The solution of Burgers equation corresponding to \eqref{BurgerIC1} remains smooth until pre-shock time $T =Tb$. In Table \ref{tabBurgersIC1},	the $L^\infty$ and $L^1$ convergence rate of different non-oscillatory entropy stable schemes is given and compared with the corresponding non-oscillatory scheme with respect to the initial condition \eqref{LinearIC1}. Similar to linear case, entropy stable schemes {\it EC-m-$F^w$-n} maintains the $k^{th}$ order convergence rate, where $k=min(m,n)$.
\item {IC4: To test non-oscillatory property:}
\begin{equation}\label{BurgerIC2}
u(x,0) = \left\{\begin{array}{cc}
\sin(\pi x), & |x|> 4,\\
3, & -1\leq x\leq -0.5,\\
1.0, & -0.5\leq x\leq 0,\\
3.0, & 0.0\leq x < 0.5,\\
1.0, & \mbox{else}.
\end{array}\right. x\in [-4,4].	
\end{equation}
The initial condition \eqref{BurgerIC2} contains smooth and discontinuous data regions which eventually develops a complex solution containing stationary shock at $x=\pm 3$, rarefaction waves and moving shocks of varied intensities. In Figure \ref{fig:burgertest6eceno}, numerical solution computed by different high order entropy stable schemes is given and compared. It can be observed from results in sub-figures \ref{fig:burgertest6eceno}(a) and \ref{fig:burgertest6eceno}(b) that entropy stable non-oscillatory {\it TeCNO} schemes proposed in  \cite{Fjordholm2012} and low dissipative {\it TeC-WENOJS3} scheme in \cite{BbRk} respectively exhibits oscillations at discontinuities. However, the proposed  entropy stable schemes {\it EC-m-$F^w$-n} in this work give non-oscillatory solution as shown in sub-figures \ref{fig:burgertest6eceno}(c) and \ref{fig:burgertest6eceno}(d). Moreover, the resolution of discontinuities is characterized by the base non-oscillatory flux used in the construction of entropy stable scheme. 
\end{itemize}

\begin{table}[htb!]\label{tabBurgersIC1}
	\centering	
	\begin{tabular}{|cc|}
		\hline
			{\bf ENO-3} & EC4-ENO-3\\
		\begin{tabular}{|c|c|c|c|c|}
		\hline N  & $L^\infty$ error  &  Rate &    $L^1$ error    & Rate \\
		\hline   20   &   4.0898e-03   &  -   &   2.5767e-03   &  -   \\
		\hline   40   &   6.5296e-04   &   2.65   &   3.3011e-04   &   2.96   \\
		\hline   80   &   1.3560e-04   &   2.27   &   4.4561e-05   &   2.89   \\
		\hline   1.60   &   2.9932e-05   &   2.18   &   6.2369e-06   &   2.84   \\
		\hline   3.20   &   6.3461e-06   &   2.24   &   9.5137e-07   &   2.71   \\
		\hline 
	\end{tabular}&	
	\begin{tabular}{|c|c|c|c|c|}
		\hline N  & $L^\infty$ error  &  Rate &    $L^1$ error    & Rate \\
		\hline   20   &   3.9841e-03   &  -   &   2.5641e-03   &  -   \\
		\hline   40   &   6.9710e-04   &   2.51   &   3.3276e-04   &   2.95   \\
		\hline   80   &   1.4126e-04   &   2.30   &   4.4795e-05   &   2.89   \\
		\hline   160   &   3.0229e-05   &   2.22   &   6.2516e-06   &   2.84   \\
		\hline   320   &   6.3530e-06   &   2.25   &   9.5370e-07   &   2.71   \\
		\hline 
	\end{tabular}\\
	 & \\
	{\bf WENOJS-3} & 	EC4-WENOJS-3\\
	\begin{tabular}{|c|c|c|c|c|}
		\hline N  & $L^\infty$ error  &  Rate &    $L^1$ error    & Rate \\
		\hline   20   &   3.0054e-02   &  -   &   2.1504e-02   &  -   \\
		\hline   40   &   1.1646e-02   &   1.37   &   5.7054e-03   &   1.91   \\
		\hline   80   &   4.4466e-03   &   1.39   &   1.2723e-03   &   2.16   \\
		\hline   160   &   1.2298e-03   &   1.85   &   2.0970e-04   &   2.60   \\
		\hline   320   &   1.6020e-04   &   2.94   &   1.9564e-05   &   3.42   \\
		\hline 
	\end{tabular}&
	\begin{tabular}{|c|c|c|c|c|}
		\hline N  & $L^\infty$ error  &  Rate &    $L^1$ error    & Rate \\
		\hline   20   &   2.8160e-02   &  -   &   2.1200e-02   &  -   \\
		\hline   40   &   1.1420e-02   &   1.30   &   5.6211e-03   &   1.92   \\
		\hline   80   &   4.3542e-03   &   1.39   &   1.2561e-03   &   2.16   \\
		\hline   160   &   1.1957e-03   &   1.86   &   2.0148e-04   &   2.64   \\
		\hline   320   &   1.8016e-04   &   2.73   &   1.6295e-05   &   3.63   \\
	\hline
	\end{tabular}\\
	 &\\	
		{\bf WENOJS-5}& 		EC6-WENOJS-5\\ 
	\begin{tabular}{|c|c|c|c|c|}
		\hline N  & $L^\infty$ error  &  Rate &    $L^1$ error    & Rate \\
		\hline   20   &   1.5856e-03   &  -   &   8.6897e-04   &  -   \\
		\hline   40   &   4.0537e-05   &   5.29   &   2.5151e-05   &   5.11   \\
		\hline   80   &   1.9176e-06   &   4.40   &   7.5007e-07   &   5.07   \\
		\hline   160   &   1.5094e-07   &   3.67   &   2.3675e-08   &   4.99   \\
		\hline   320   &   3.9001e-09   &   5.27   &   5.8722e-10   &   5.33   \\
		\hline 
	\end{tabular}&
	\begin{tabular}{|c|c|c|c|c|}
		\hline N  & $L^\infty$ error  &  Rate &    $L^1$ error    & Rate \\
		\hline   20   &   1.5038e-03   &  -   &   8.3906e-04   &  -   \\
		\hline   40   &   4.0904e-05   &   5.20   &   2.4423e-05   &   5.10   \\
		\hline   80   &   1.9740e-06   &   4.37   &   7.3608e-07   &   5.05   \\
		\hline   160   &   1.5348e-07   &   3.68   &   2.3303e-08   &   4.98   \\
		\hline   320   &   4.0946e-09   &   5.23   &   5.7234e-10   &   5.35   \\
		\hline		
	\end{tabular} \\ & \\
\hline
\end{tabular}
	\caption{{\it Burgers equation:} Convergence rate of base non-oscillatory schemes and corresponding entropy stable schemes for initial condition \eqref{BurgerIC1}, $CFL=0.4,\;\Delta t = \Delta x^{\frac{5}{3}},\;  T_f =\frac{1}{2\pi}$. }
\end{table} 
   
\begin{figure}
	\begin{tabular}{cc}	
		\includegraphics[scale=0.5]{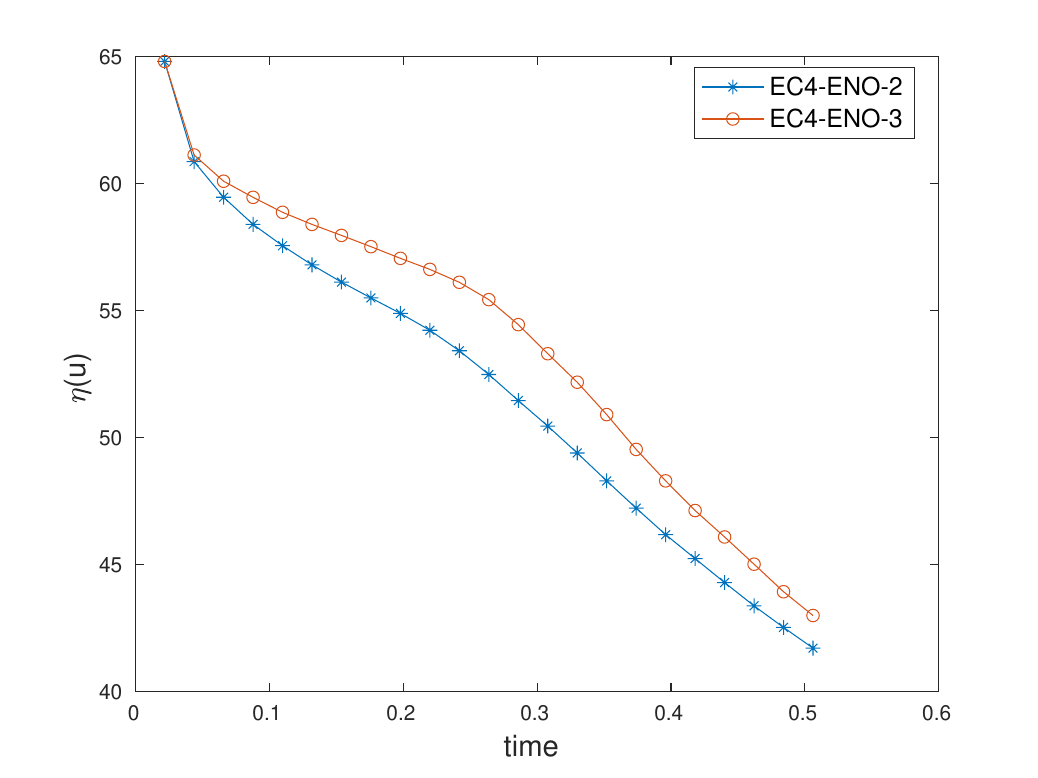} & \includegraphics[scale=0.5]{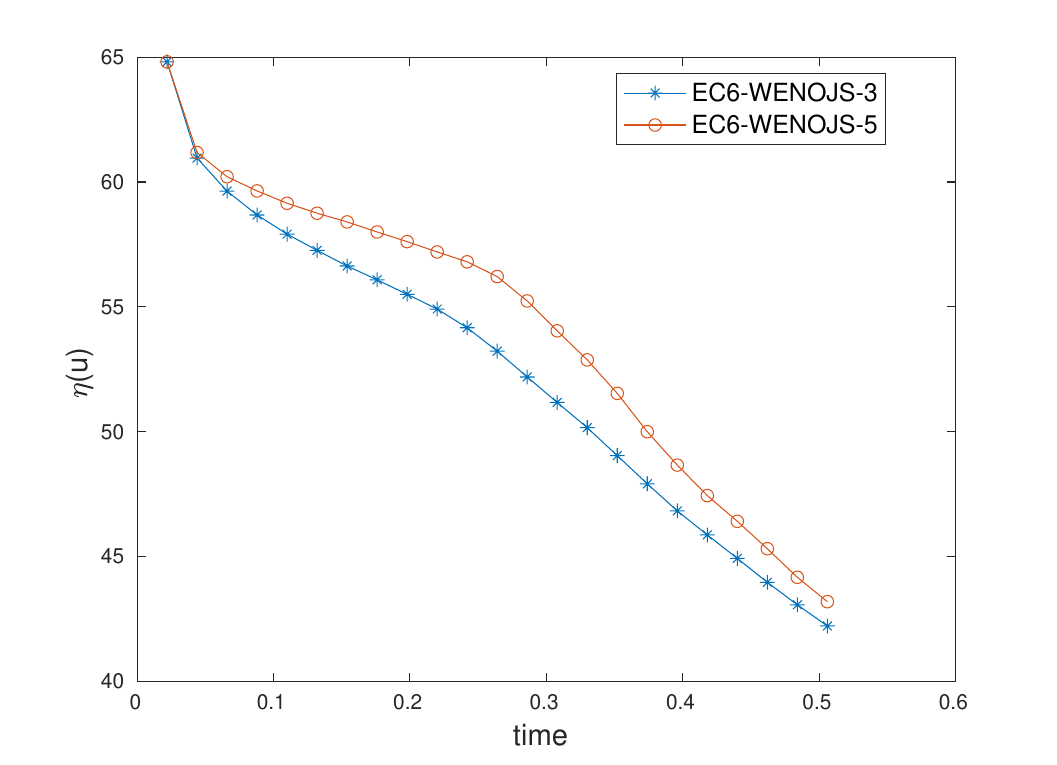}\\
		(a) & (b)
	\end{tabular}
	\caption{Entropy decay corresponding to the computed solution in Figure \ref{fig:burgertest6eceno} by {\it EC-m-$F^w$-n} schemes.}
	\label{fig:burgertest6entropy}
\end{figure}

\newpage
\subsection{Systems of Conservation Laws}
Like scalar case, in system case also various proposed entropy stable schemes {\it EC-m-$F^w$-n}  are used to compute the numerical solution of hyperbolic systems. However to restrict the length of the presentation, in one dimensional tests, results are reported only for third and fifth order entropy stable {\it ES-4-ENO-3}  and {\it ES-6-WENOJS-5} schemes. Also results are compared with the result by base non-oscillatory {\it ENO-3} and {\it WENOJS-5} schemes. In order to compare the dissipation and discontinuities resolution nature of the proposed schemes, comparison of numerical results by {\it EC-6-WENOJS-5} is done  with {\it TECNO5} flux \cite{Fjordholm2012} which is constructed from the sum of an EC-6 flux \eqref{ECflux6} and an entropy dissipative term using local Lax-Friedrichs dissipation the entropy stable Roe-type dissipation \cite{Chandrashekar2013}. These results are given in \ref{fig:eulertest1}(c), \ref{fig:eulertest1}(d) and \ref{fig:eulertest3}(c), \ref{fig:eulertest3}(d). In two dimensional test cases, computational results by fifth order non-oscillatory entropy stable scheme {\it EC6-WENOJS-5} are given and compared with the result of non-oscillatory WENO  {\it WENOJS-5} schemes.
\subsubsection{The 1D Euler system}  The one dimensional system of Euler equations is given by 
\begin{equation}
\left(
\begin{array}{c}
	\rho\\
	\rho u\\
	E
\end{array}
\right)_t+\left(
\begin{array}{c}
	\rho u\\
	\rho u^2+p\\
	u(E+p)
\end{array}
\right)_x=0,
\end{equation}
where following relationship holds between density ($\rho$), pressure($p$) and energy($E$)
\begin{equation}
	p=(\gamma-1)\bigg(E-\frac{1}{2}\rho u^2\bigg),	
\end{equation}
where $\gamma$ is the ratio of specific heat coefficient. Computational results are obtained for various Riemann problem of the form \cite{toro2013riemann}
\begin{equation}
\mathbf{u}(x,0)=\left\{\begin{array}{cl}
\mathbf{u}_l, & \text{if} \;\; x<x_0, \\
\mathbf{u}_r, & \text{if} \;\; x\geq x_0,
\end{array}
		\right.
\end{equation}
where $\mathbf{u}_l=(\rho_l,u_l,p_l)$ and $\mathbf{u}_r=(\rho_r,u_r,p_r)$.\\\\
\subsubsection*{Sod shock tube test:} 
\begin{figure}
\hspace{-0.5cm}
	\begin{tabular}{cc}\includegraphics[scale=0.55]{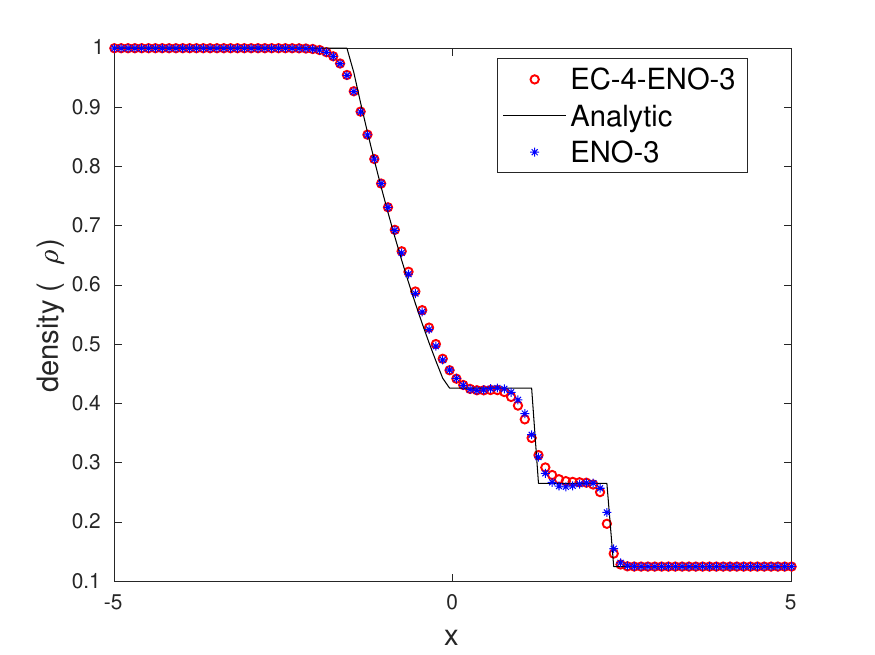} &
	\includegraphics[scale=0.55]{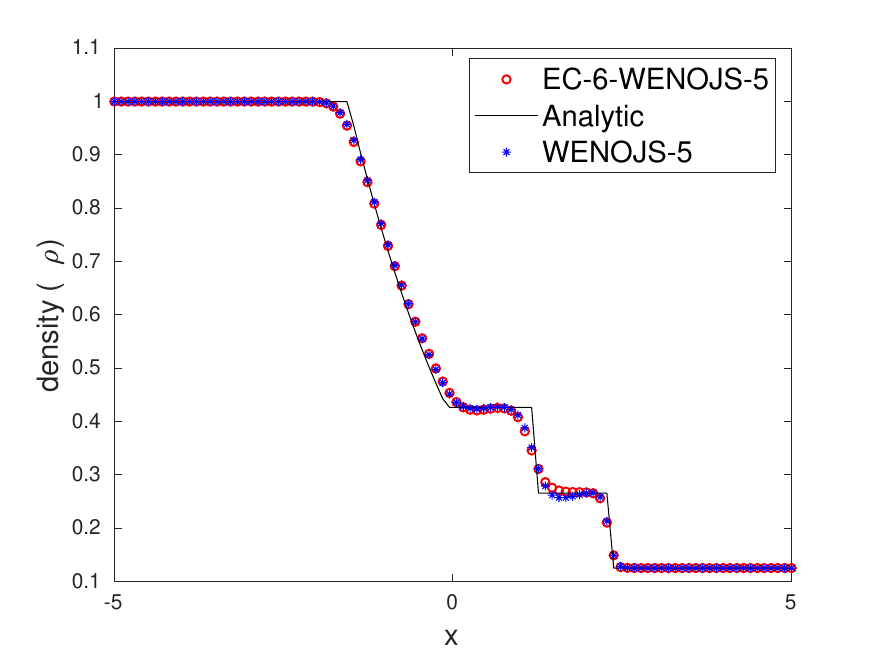}\\
	(a) & (b)\\ 
	\includegraphics[scale=0.45]{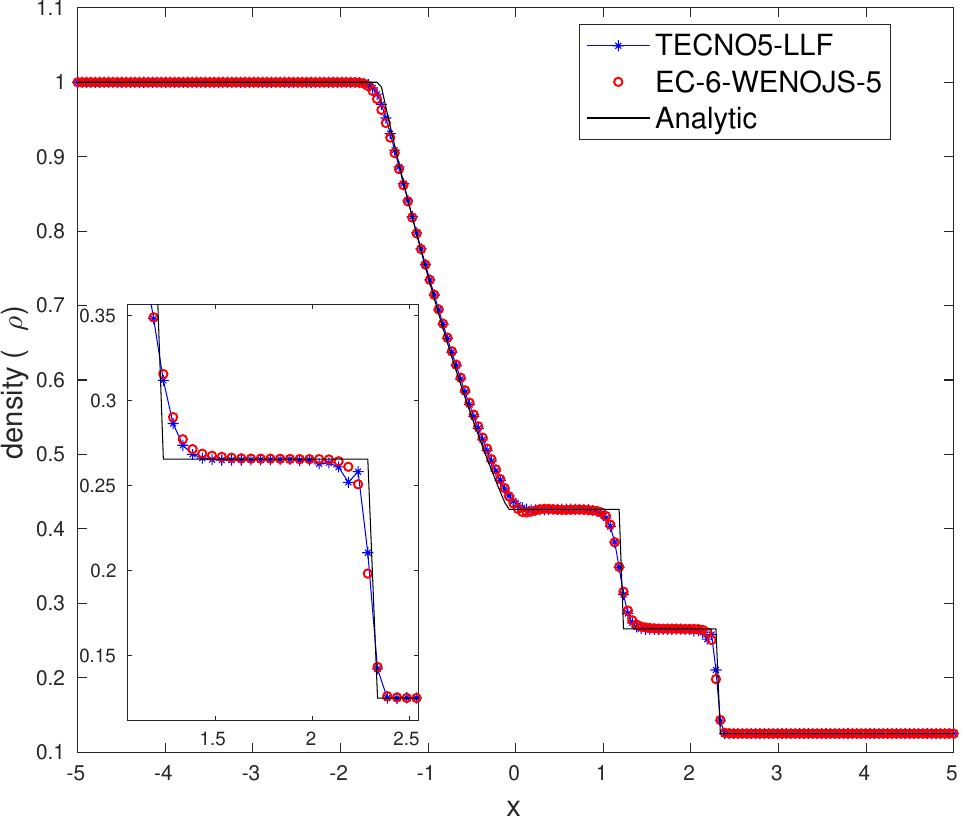} &
	\includegraphics[scale=0.45]{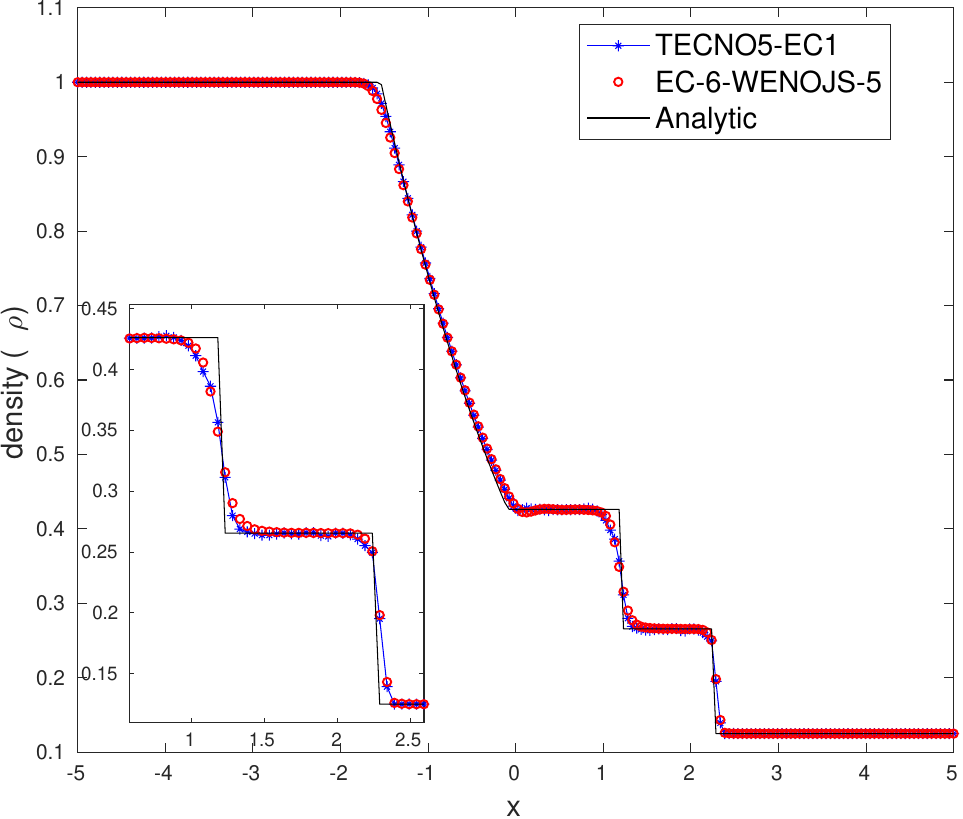}\\	
	(c) & (d)\\
	\includegraphics[scale=0.55]{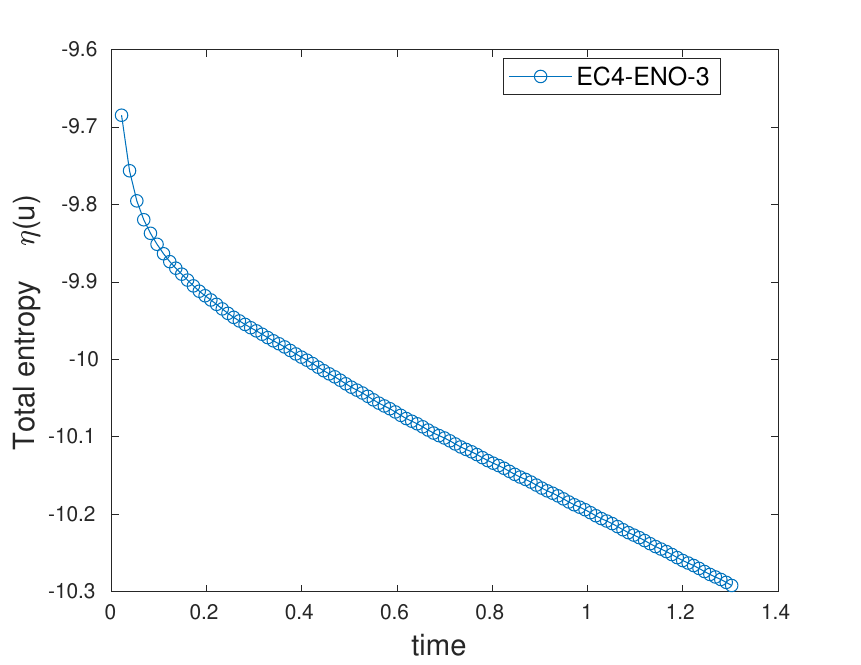} &
	\includegraphics[scale=0.55]{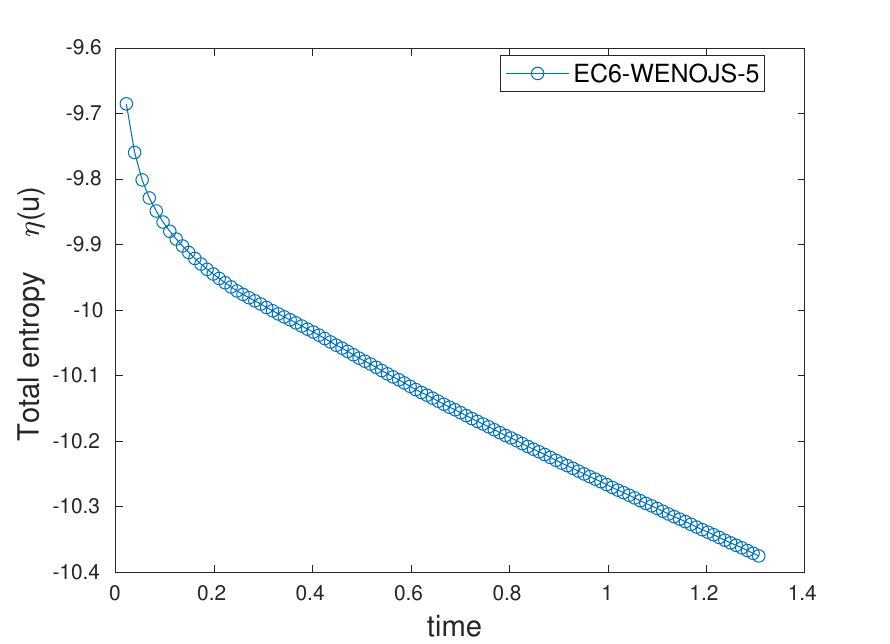}\\
	(e) & (f)
	\end{tabular}
	\caption{Sod Shock tube solution (a) \& (b) at $N=100,\; CFL=0.25,\; T_f=1.3$. Comparison of dissipative nature of {\it EC-6-WENOJS-5} with {\it TECNO5} scheme (c) and (d)}. Entropy $\eta(\mathbf{u})$ decay (e) \& (f) respectively where $\eta(\mathbf{u}) = \displaystyle \frac{\rho(\log(p)-\gamma \log(\rho))}{\gamma -1}$ as defined in \cite{Ismail2009}. 
	\label{fig:eulertest1}
\end{figure}

The sod problem is defined in \cite{sod1978survey} is given by following initial condition
\begin{equation}
	\mathbf{u}(x,0)=\left\{\begin{array}{cl}
	(1,0,1), & \text{if} \;\; -5<x<0, \\
	(0.125,0,0.1), & \text{if} \;\; 0 \leq x \leq5.
	\end{array}
	\right.
\end{equation}
For this initial condition the evolved solution consists of a rarefaction wave followed by a contact discontinuity and the shock discontinuity. The numerical solution for this test is given and compared in Figure \ref{fig:eulertest1}. 

%

\subsubsection*{Laney Shock tube Test} 
\begin{figure}
	\hspace{-1cm}
	\begin{tabular}{cc}\includegraphics[scale=0.65]{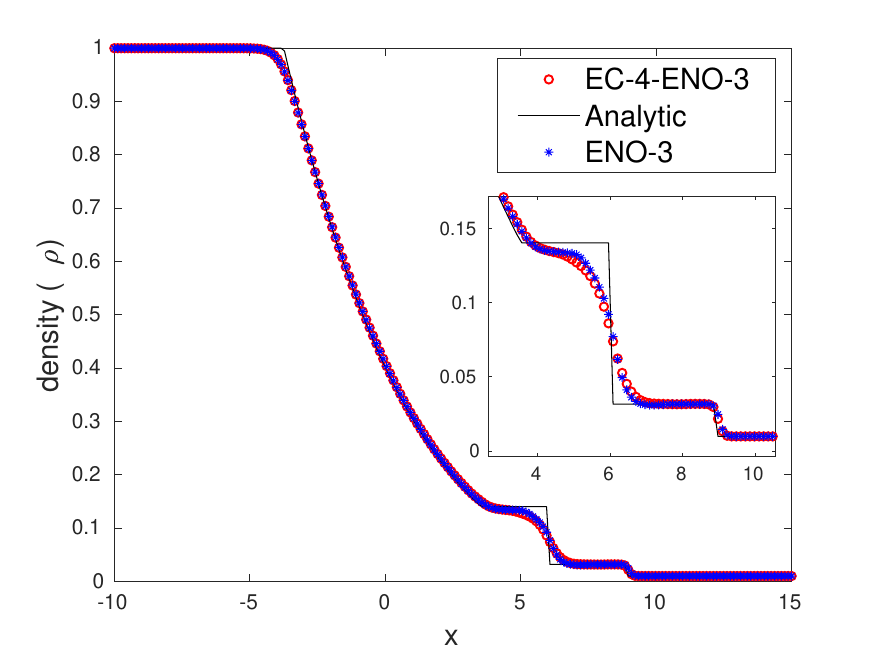}&\includegraphics[scale=0.65]{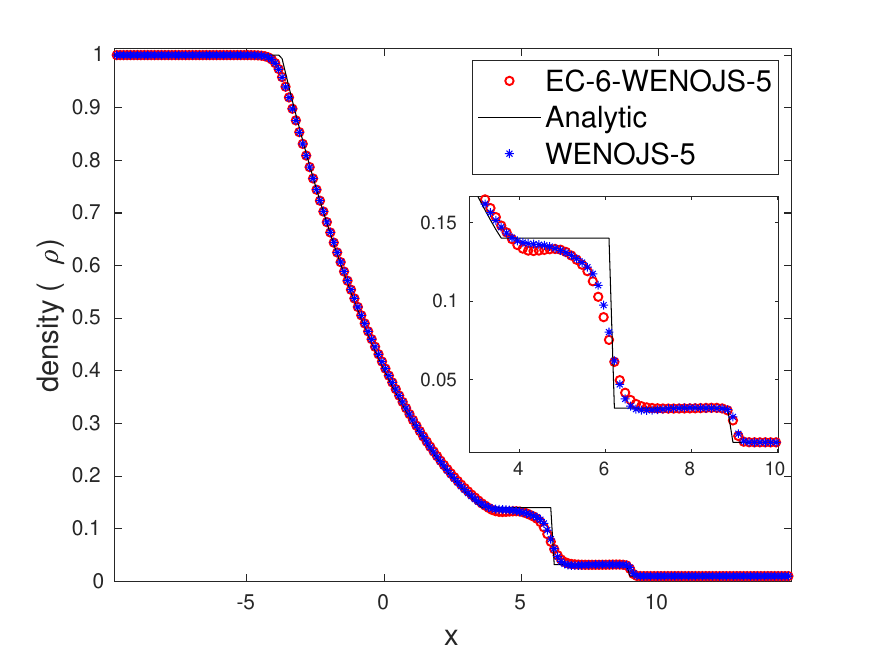}
	\end{tabular}
	\caption{Density plot for Laney test with $N=200, CFL=0.25, T_f= 0.01$.}
	\label{fig:eulertest2}
\end{figure}

A more complicated shock tube test problem to check the non-oscillatory property of any numerical scheme is taken from \cite{Laney} which is also known as the Leblanc shock tube and govern by the initial Riemann data  
\begin{equation}
	\mathbf{u}(x,0)=\left\{\begin{array}{cl}
	(1.0,0,100000), & \text{if} \;\; -10\leq x<0, \\
	(0.01,0,1000), & \text{if} \;\; 0 \leq x \leq 10.
	\end{array}
	\right.
\end{equation}
In this test, the ratio of left and right states of density and pressure across initial discontinuity is very high and right initial state of density is close to zero. Therefore, computationally, even small oscillations can lead to negative density or pressure, which results into nonphysical imaginary speed of sound $c= \sqrt{\frac{\gamma p}{\rho}}$. The numerical solution for this test case is given in Figure \ref{fig:eulertest2}.
\subsubsection*{Lax tube test:}
\begin{figure}
	\hspace{-1cm}
	\begin{tabular}{cc}\includegraphics[scale=0.65]{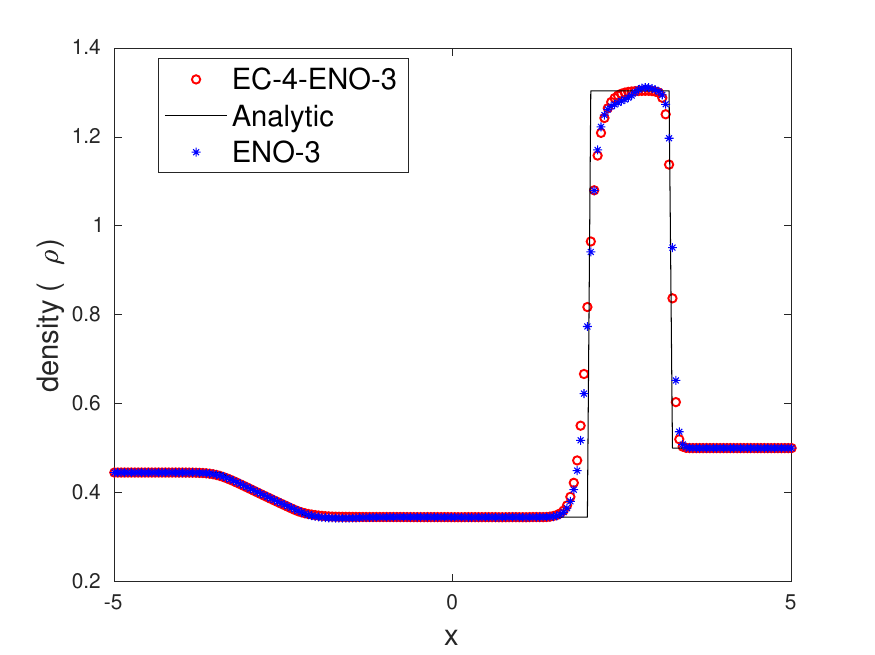}& \includegraphics[scale=0.65]{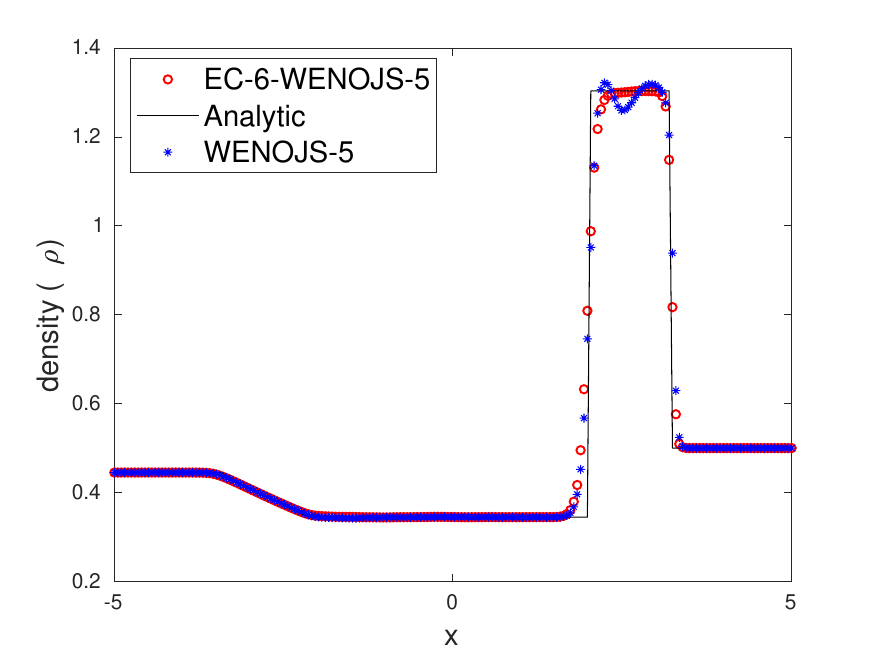}\\
		(a) & (b)\\
		\includegraphics[scale=0.45]{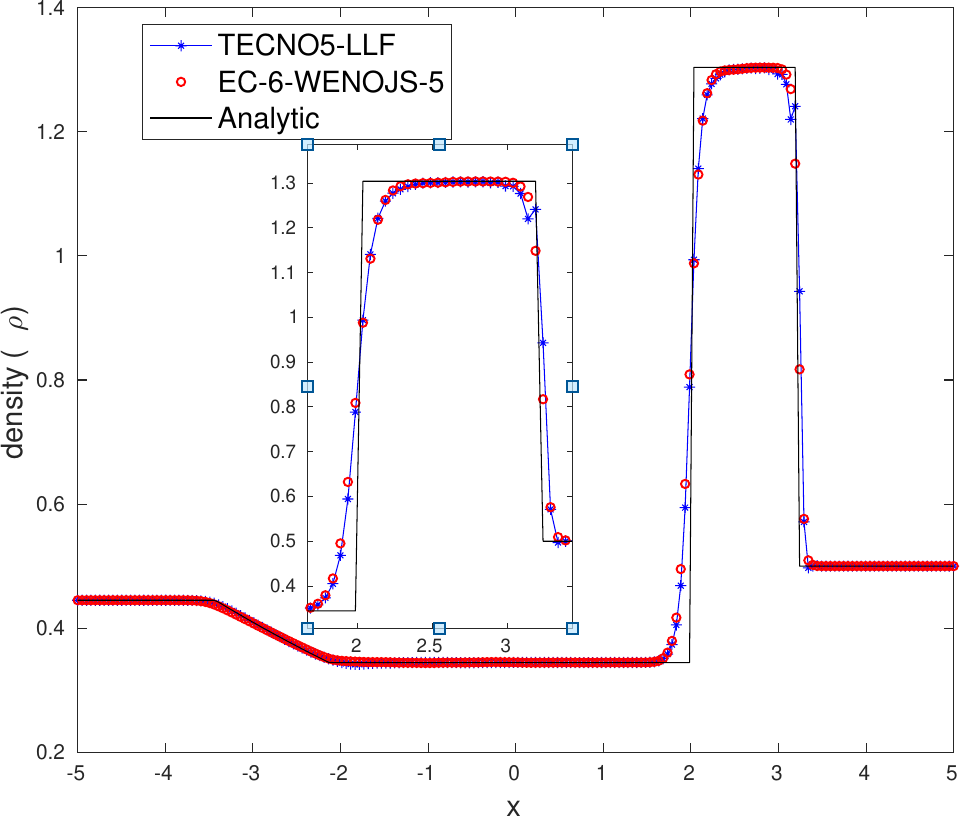} &
		\includegraphics[scale=0.45]{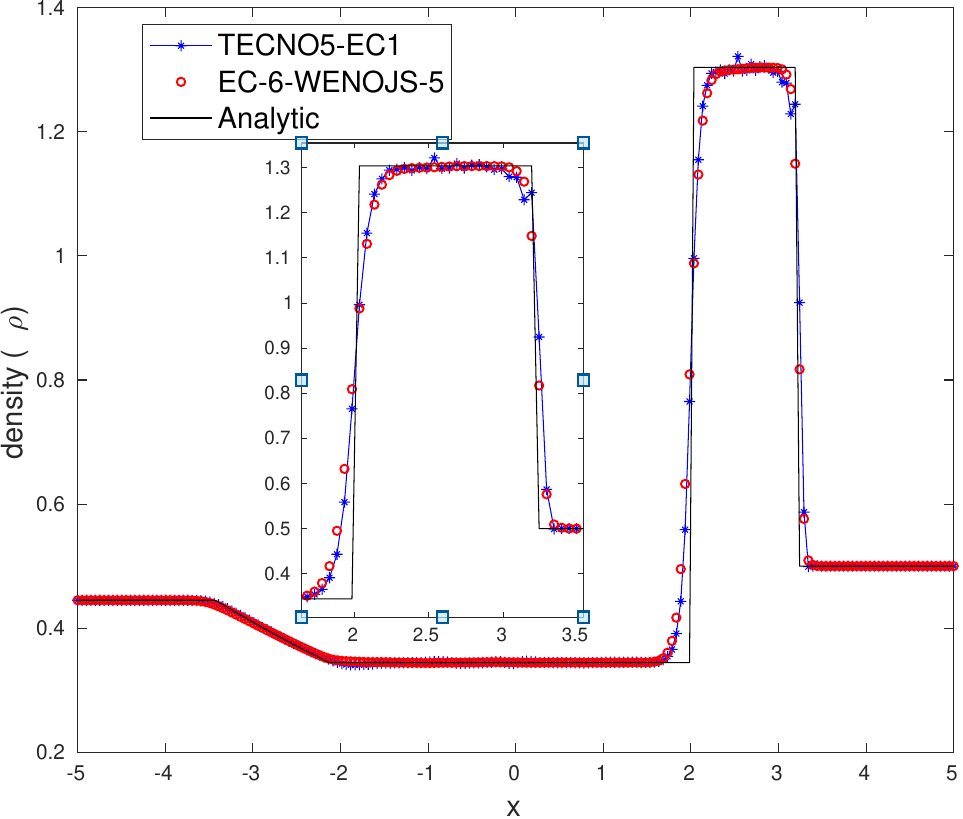}\\	
		(c) & (d)\\
	\end{tabular}
	\caption{Solution at $N=200, CFL=0.25, T_f=1.3$. Comparison of discontinuity resolving nature of {\it EC-6-WENOJS-5} with {\it TECNO5} scheme (c) and (d).}
	\label{fig:eulertest3}
\end{figure}

We consider the Lax tube problem discussed in \cite{Lax1954} with the initial condition given by
\begin{equation}
	\mathbf{u}(x,0)=\left\{\begin{array}{cl}
	(0.445,0.698,3.528), & \text{if} \;\; -5<x<0, \\
	(0.5,0,0.571), & \text{if} \;\; 0 \leq x \leq 5.
	\end{array}
	\right.
\end{equation}
Solution corresponding this initial condition contains a right traveling strong shock wave, a contact surface, and a left rarefaction wave. The numerical solutions for this test case is given in Figure \ref{fig:eulertest3} which clearly show that the underlying entropy stable scheme completely removes the small oscillations exhibited by non-oscillatory scheme WENOJS-5 in this test case.

\subsubsection*{Arora and Roe test}
\begin{figure}[H]
	\hspace{-1cm}	
	\begin{tabular}{cc}\includegraphics[scale=0.65]{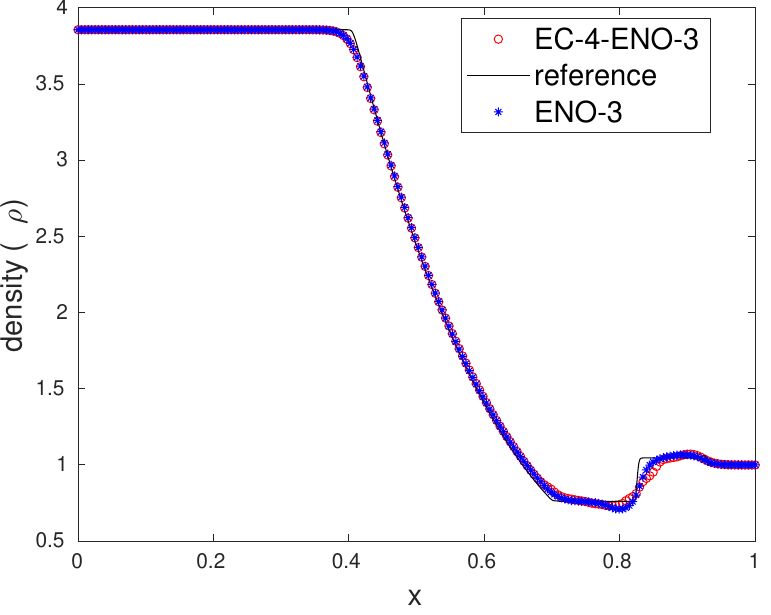}& 
		\includegraphics[scale=0.65]{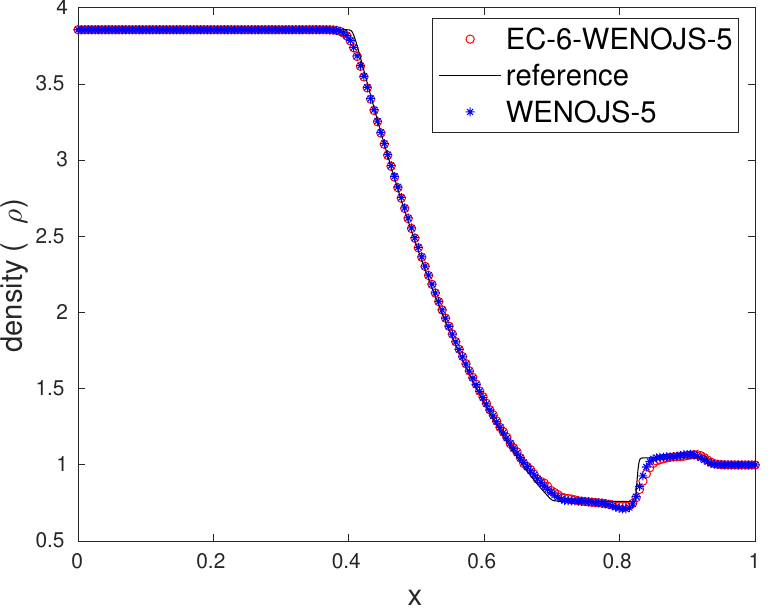}
	\end{tabular}
	\caption{\textcolor{blue}{Density plot for Mach 3 test with strong expansion fan, $N=200, CFL = 0.25, T_f=0.09$}}
	\label{fig:eulertest10}
\end{figure}

A Mach 3 test case where supersonic flow occurs is considered by Arora and Roe in \cite{Arora1997} which corresponds to the following initial condition.
\begin{equation}
	\mathbf{u}(x,0)=\left\{\begin{array}{cl}
		(3.857,0.92, 10.333), & \text{if} \;\; 0<x<0.5, \\
		(1.0,3.55,1.0), & \text{if} \;\; 0.5 \leq x \leq 1.
	\end{array}
	\right.
\end{equation}
The flow in this test is dominated by a strong expansion fan \cite{wessling}. The numerical solution for this test case is given and compared with reference solution in Figure \ref{fig:eulertest10}.
\subsubsection*{\bf Shock-entropy wave interaction} 
\begin{figure}[H]
	\hspace{-1cm}
	\begin{tabular}{cc}\includegraphics[scale=0.65]{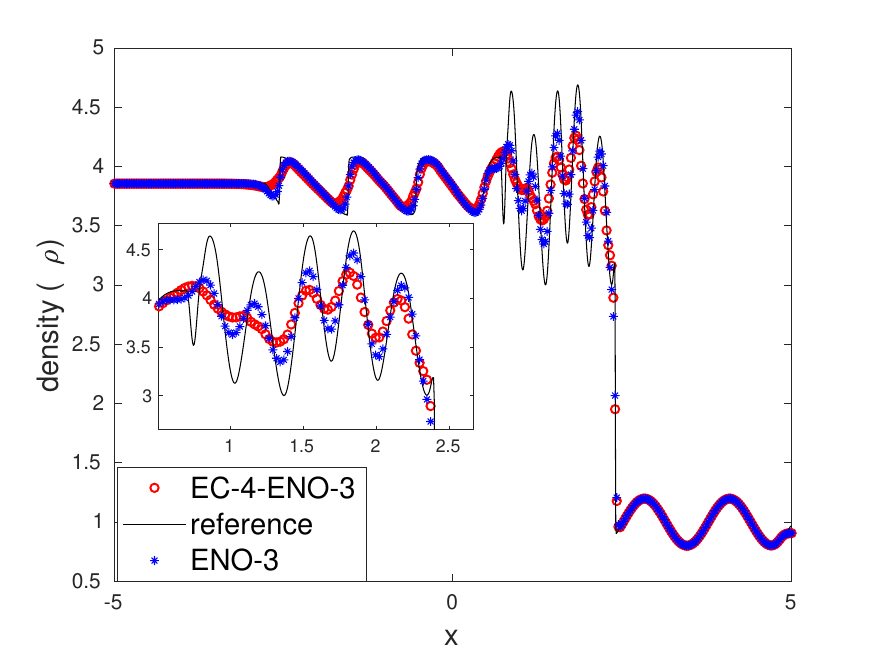}& 
		\includegraphics[scale=0.65]{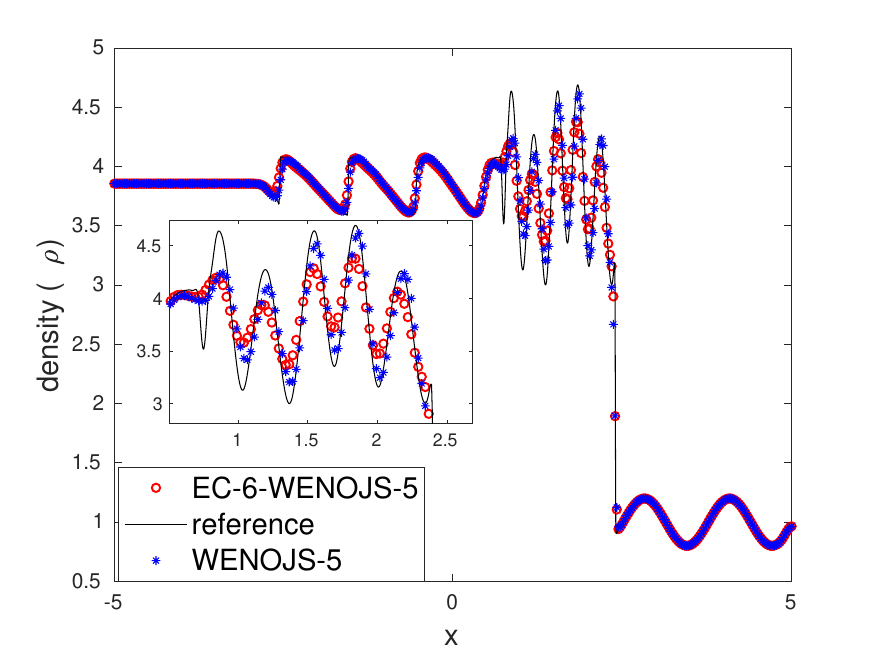}
	\end{tabular}
	\caption{Density plot for Shu-Osher test with $N=400, CFL=0.25, T_f= 1.8$. }
	\label{fig:eulertest4}
\end{figure}

This Shu-Osher problem in \cite{SHU1988439} is governed by the following initial condition 
\begin{equation}
		\mathbf{u}(x,0)=\left\{\begin{array}{cl}
			(3.857143,2.629369,2.629369), & \text{if} \;\; -5<x<0, \\
			(1+\epsilon sin(kx),0,1), & \text{if} \;\; 0 \leq x \leq5.
		\end{array}
		\right. \epsilon =0.2, k=5,
\end{equation}
which simulates shock-turbulence interaction in which a strong shock wave propagates into density field with artificial fluctuations with amplitude $\epsilon=0.2$ and wave number $k=5$. This problem tests the capability of any scheme to accurately capture a shock wave, its interaction with an unsteady density field, and the sinusoidal waves propagating downstream of the shock. The numerical solution for this test case is given in Figure \ref{fig:eulertest4}.
\subsubsection*{Two interacting blast wave}
\begin{figure}[H]
\hspace{-1cm}	
\begin{tabular}{cc}\includegraphics[scale=0.65]{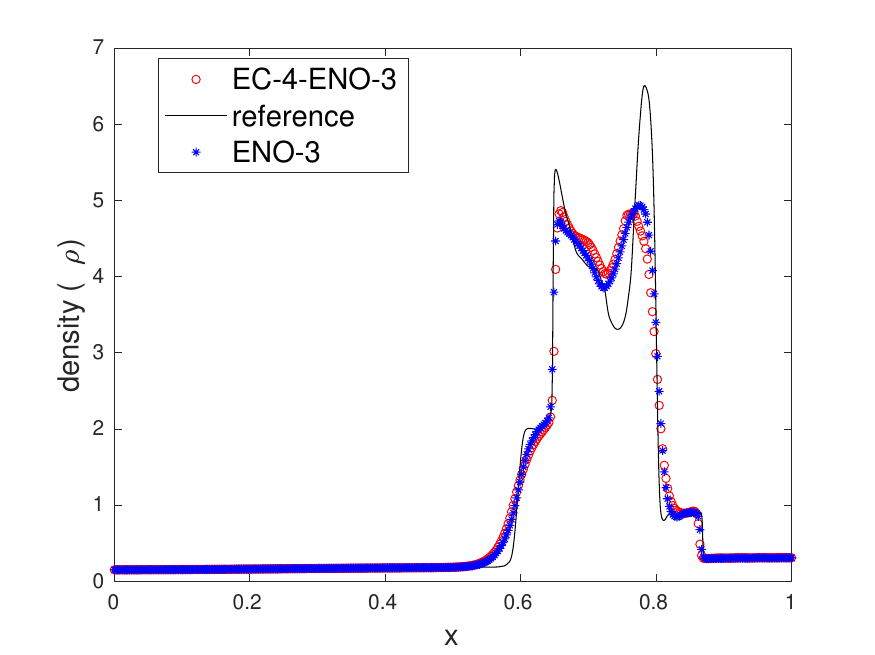}& 
	\includegraphics[scale=0.65]{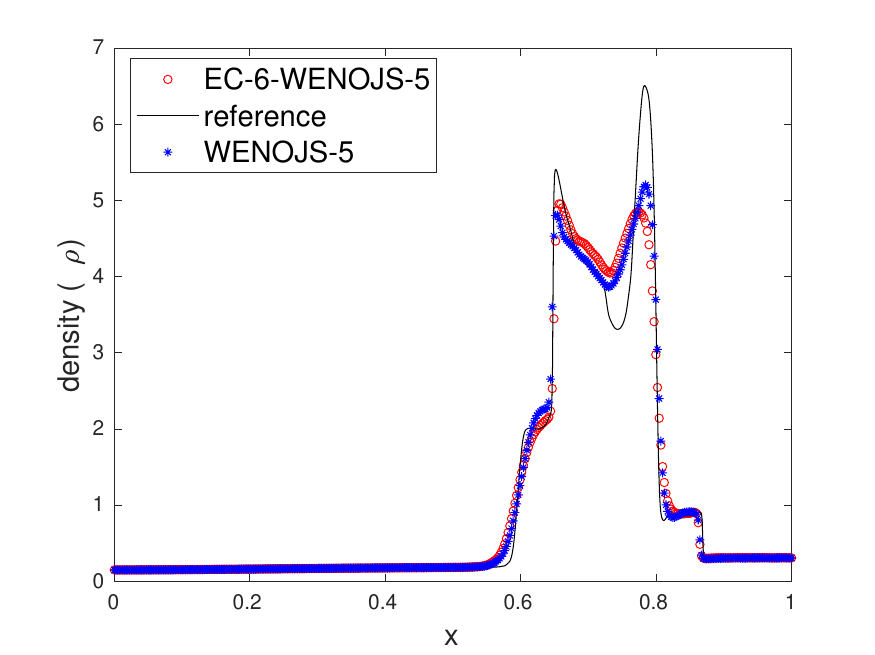}
\end{tabular}
	\caption{Density plot for two interacting blast wave test $N=400, CFL = 0.25, T_f=0.038$}
	\label{fig:eulertest9}
\end{figure}

Woodward-Colella blast wave \cite{woodward1984numerical} is another interesting  problem to test the shock capturing ability of numerical scheme given by the following initial condition.
\begin{equation}
		\mathbf{u}(x,0)=\left\{\begin{array}{cl}
			(1.0,0.0,1000.0), & \text{if} \;\; 0.0<x<0.1, \\
			(1.0,0.0,0.01) & \text{if} \;\; 0.1 \leq x \leq 0.9, \\
			(1.0,0.0,100.0), & \text{if} \;\; 0.9 \leq x \leq 1.0 .
		\end{array}
		\right.
\end{equation}
This problem involves the multiple interactions of shock, contact, and rarefaction wave. A  reflecting boundary condition is applied at the boundaries $x=0$ and $x=1$ of the domain. The numerical solution for this test case is given in Figure \ref{fig:eulertest9}.
It can be concluded from the numerical results for 1D shock tube problems in figures \ref{fig:eulertest1}-\ref{fig:eulertest9}, that high order entropy stable schemes {\it EC-4-ENO-3} and {\it EC-6-WENOJS-5} yield results comparable with essentially non-oscillatory ENO/WENO schemes even for complex flow with strong shocks and expansion waves. Also results in Figures  \ref{fig:eulertest1}(c), \ref{fig:eulertest1}(d) and \ref{fig:eulertest3}(c), \ref{fig:eulertest3}(d) show that the dissipation and shock resolution of the proposed scheme {\it EC-6-WENOJS-5} is comparable with {\it TECNO5} schemes.
\subsection{2D Euler equation} 
 In this section the proposed scheme is applied to 2D Euler equations
        \begin{equation}\label{Euler2Deq}
        	\left(
        	\begin{array}{c}
        		\rho\\
        		\rho u\\
        		\rho v\\
        		E
        	\end{array}
        	\right)_t+\left(
        	\begin{array}{c}
        		\rho u\\
        		\rho u^2+p\\
        		\rho uv \\
        		u(E+p)
        	\end{array}
        	\right)_x+\left(
        	\begin{array}{c}
        		\rho v\\
        		\rho uv \\
        		\rho v^2+p\\
        		
        		v(E+p)
        	\end{array}
        	\right)_y=0,
        \end{equation}
    where $\rho$ is density and $u,\,v$ are component of velocity along $x$ and $y$ direction respectively.The pressure and energy are related by the following 
\begin{equation} 	
E=\frac{p}{\gamma-1}+\frac{\rho(u^2+v^2)}{2},
\end{equation}
and $\gamma$ is the ratio of specific heat. 
The following test problems are considered for testing the performance of the scheme for equation \eqref{Euler2Deq}. In all the tests numerical solution using fifth order non-oscillatory entropy stable scheme {\it EC6-WENOJS-5} is compared with the base non-oscillatory schemes {\it WENOJS-5}. Solution plots are given in the figure side by side.
\par
\begin{figure}

\begin{tabular}{cc}
	\vspace{-8cm}\\ 
	\hspace{-3cm}\includegraphics[width=0.8\linewidth]{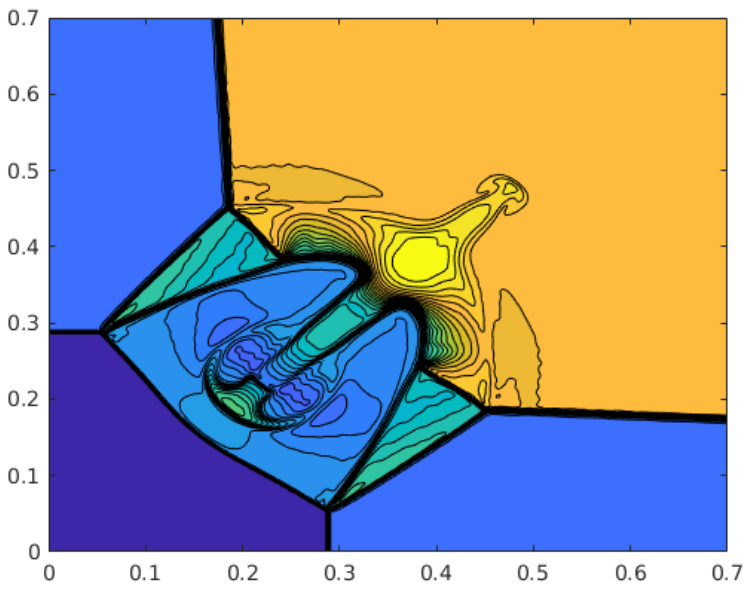} & 
	\hspace{-5cm}\includegraphics[width=0.8\linewidth]{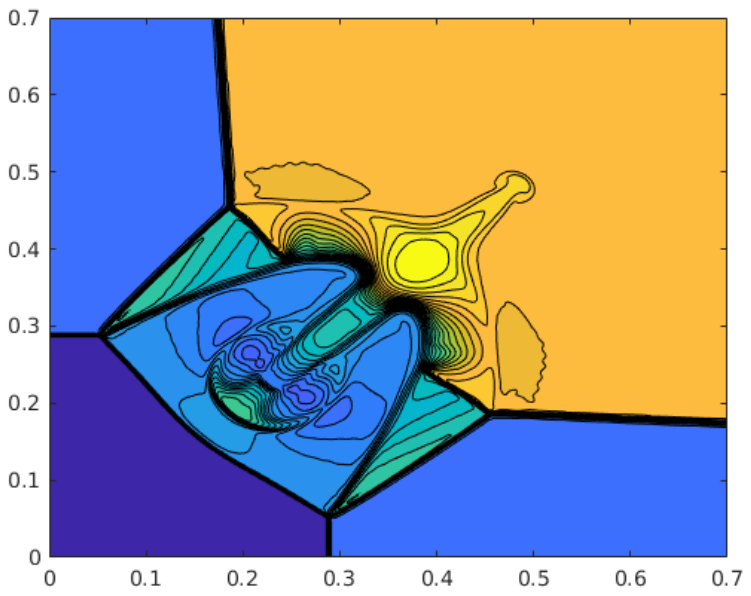}\vspace{-6cm}\\
	\hspace{-2cm}(a) & \hspace{-5cm}(b)\\ 
\end{tabular}
\caption{ Contour solution plot with $31$ of contour lines corresponding to configuration \eqref{config3} by {\it WENOJS-5} in (a) and by {\it EC6-WENOJS-5} in (b), $CFL=0.25,\; T_f=0.5, N=400\times 400$.   }
\label{fig:euler2dtest3}
\end{figure}

 {\bf 2D Riemann problem \cite{Kurganov2002, schulz1993numerical}} Consider 2D Euler equations \eqref{Euler2Deq} with Riemann data defined in 
in the following way,
 \begin{eqnarray}\label{config3}
 		(p,\rho,u,v)=\left\{\begin{array}{cl}
 			(1.5000,1.5000,0.0000,0.0000), & \text{if} \;\; x>0.5 \;\; \text{and} \;\; y>0.5,\\
 			(0.3000,0.5323,1.2060,0.0000), & \text{if} \;\;  x <0.9 \;\;  \text{and} \;\;y\geq0.5, \\
 			(0.0290,0.1380,1.2060,1.2060), & \text{if} \;\; x <0.5 \;\;  \text{and} y<0.5, \\
 			(0.3000,0.5323,0.0000,1.2060), & \text{if} \;\; x>0.5 \;\; \text{and} \;\; y<0.5.\\
 		\end{array}
 		\right.
 		\end{eqnarray}
Solution for this initial condition is computed at time $t=0.5$ and corresponding filled contour plots are given in figure \ref{fig:euler2dtest3}. 

\par
\begin{figure}[H]
	\begin{tabular}{cc}
		\vspace{-8cm}\\ 
		\hspace{-3cm}\includegraphics[width=0.8\linewidth]{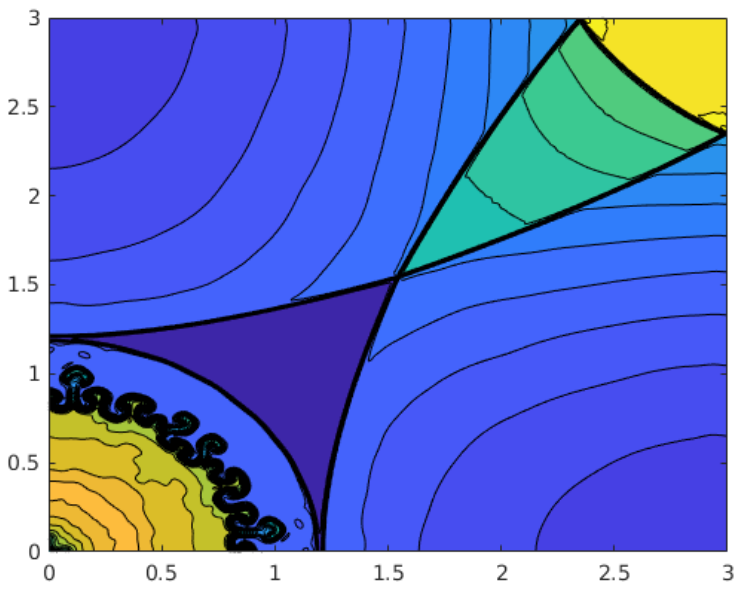} & 
		\hspace{-5cm}\includegraphics[width=0.8\linewidth]{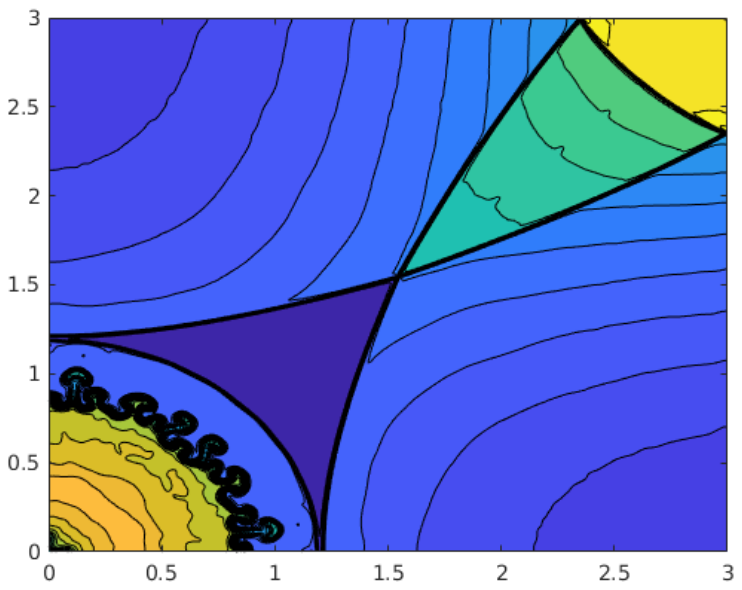}\vspace{-6cm}\\
		\hspace{-2cm}(a) & \hspace{-5cm}(b)\\ 
	\end{tabular}
\caption{Contour solution plot with $31$ contour lines corresponding to explosion problem \eqref{ICexplosion}. Solution by {\it WENOJS-5} in (a) and by {\it EC6-WENOJS-5} in (b), $CFL=0.45,\; T_f=3.2, N=400\times 400$.}
\label{fig:euler2dexplosion}
\end{figure}

 {\bf Explosion problem \cite{liska2003comparison}}
The explosion test problem is setup in a square domain $[-3,3]\times[-3,3]$ in x-y plane. The initial Riemann data is separated in the domain by a circle with center $(0,0)$ and radius $0.4$. The initial density and pressure are defined in the following way.
\begin{equation}\label{ICexplosion}
\left\{\begin{array}{cl}
\rho(x,y)=1,p(x,y)=1,  & \text{if} \;\; x^2+y^2<(0.4)^2,\\
\rho(x,y)=0.125,p(x,y)=0.1, \; &\text{otherwise}.
\end{array}
\right.
\end{equation}
 The filled contour plot of numerical results for explosion problem are shown and compared in figure \ref{fig:euler2dexplosion}.  
\par
\begin{figure}
\begin{tabular}{cc}
	\vspace{-8cm}\\ 	
	\hspace{-3cm}\includegraphics[width=0.8\linewidth]{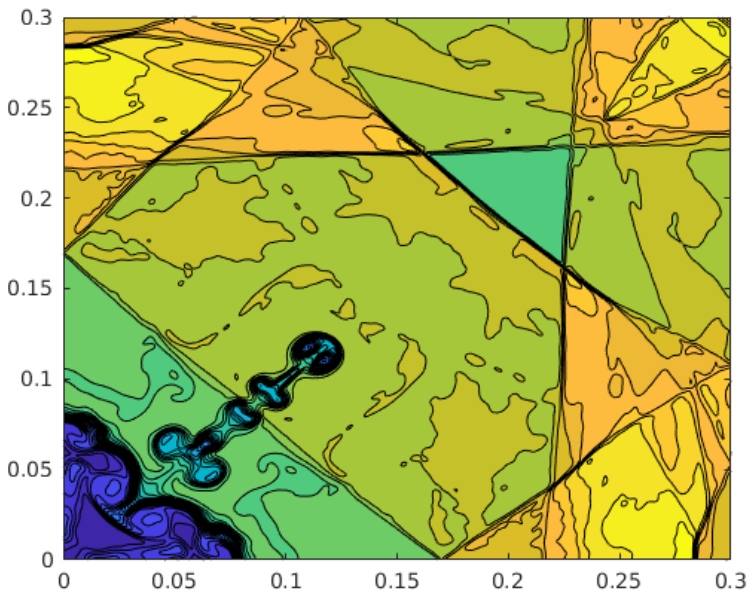} & 
	\hspace{-5cm}\includegraphics[width=0.8\linewidth]{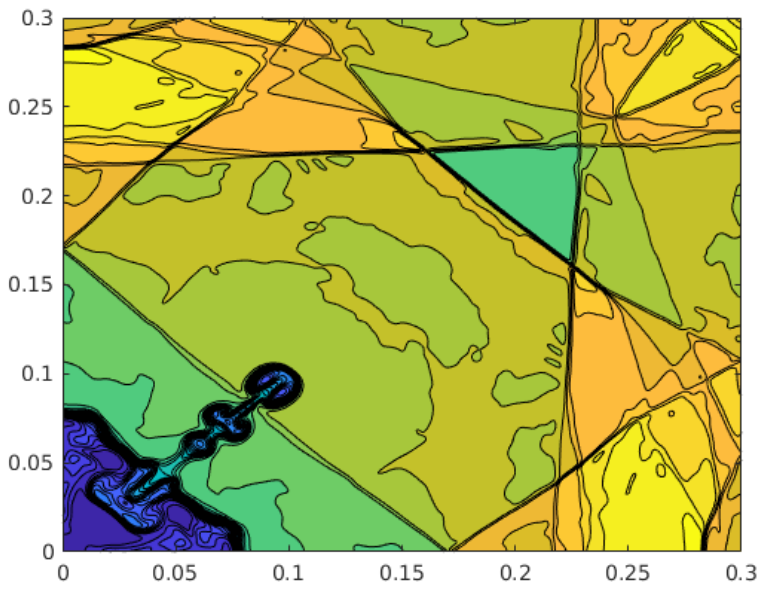}\vspace{-6cm}\\
	\hspace{-2cm}(a) & \hspace{-5cm}(b)\\ 
\end{tabular}
\caption{Contour solution plot with $31$ contour lines corresponding to implosion problem \eqref{imp}. Solution by {\it WENOJS-5} is in (a) and by {\it EC6-WENOJS-5} in (b), $CFL=0.25,\; T_f=2.5, N=400\times 400$   }
\label{fig:euler2dimplosion}
\end{figure}

 {\bf Implosion}\cite{Hui1999unified}
Consider the implosion problem modeled inside a square domain  $\left[-0.3,0.3\right]\times \left[-0.3,0.3\right]$ in $x-y$ plane. Initial Density and pressure distribution of the gas are following,
\begin{equation}\label{imp}
\begin{cases}
\rho (x,y)=0.125, p(x,y)=0.14, & \text{if $|x|+|y|<0.15$ },\\
\rho (x,y)=1, p(x,y)=1, & \text{otherwise }.
\end{cases}
\end{equation}
Initially the velocities are kept zero in the computational domain $\left[0,0.3\right]\times \left[0,0.3\right]$ with reflecting boundary. Computation is done only for the upper right quadrant $(x,y) \in (0,0.3)\times(0, 0.3)$ as in \cite{liska2003comparison, biswas2018accuracy}. The numerical results are shown and compared in Figure \ref{fig:euler2dimplosion}.
From the numerical results for 2D Euler tests in figures \ref{fig:euler2dtest3}-\ref{fig:euler2dimplosion} it is clear that {\it EC6-WENOJS-5} captures feature of the flow and results are comparable to WENO5JS scheme.
\section{Conclusion}\label{sec7}
In this work, problem of constructing non-oscillatory arbitrary order entropy stable flux is solved by framing it as least square optimization problem. Based on optimization, entropy stable flux is proposed which utilizes a {\em flux sign stability property}. Some of the existing entropy stable fluxes are retrospectively shown to satisfy this flux sign stability property. The proposed approach is robust and works well with any entropy conservative and non-oscillatory flux. Numerical results also established that such constructed entropy stable schemes give excellent non-oscillatory results even for complex problems with slightly more diffusion compared to underlying non-oscillatory scheme. Moreover, the non-oscillatory nature of the resulting entropy stable schemes is characterized by the non-oscillatory flux $F^{s}$ and for smooth solution region these schemes retain formal order of accuracy of lower order flux used in the construction.       

   \end{document}